\documentclass[a4paper,12pt]{article}
\usepackage{latexsym,amsmath,amsthm,amssymb}
\usepackage{a4wide}
\usepackage{hyperref}
\usepackage{marginnote}
\usepackage{color}
\hypersetup{
pdftitle={Improvement for Adams inequalities}   
pdfauthor={Van Hoang Nguyen},
colorlinks = true,
linkcolor = magenta,
citecolor = blue,
}

\theoremstyle{plain}
\newtheorem{theorem}{Theorem}[section]

\newtheorem*{Theorema}{Theorem A}
\newtheorem*{Theoremb}{Theorem B}
\newtheorem*{Theoremc}{Theorem C}
\newtheorem*{Theoremd}{Theorem D}

\newtheorem{proposition}[theorem]{Proposition}
\newtheorem{lemma}[theorem]{Lemma}

\theoremstyle{definition}

\theoremstyle{remark}

\renewcommand{\thefootnote}{\arabic{footnote}}
\newcommand{\C}[1]{\ensuremath{{\mathcal C}^{#1}}} %espace C^1
\def\R{\mathbb R}% tap so thuc
\def\N{\mathbb N}% tap so tu nhien
\def\C{\mathbb C}% tap so phuc
\def\sign{\mathop{\rm sign}\nolimits}%toan tu dau
\def\al{\alpha}% alpha
\def\om{\omega}% omega
\def\Om{\Omega}% Omega
\def\be{\beta}% beta
\def\de{\delta}% delta
\def\ep{\epsilon}% epsilon
\def\na{\nabla}% nabla
\def\lt{\left}% trai
\def\rt{\right}% phai
\def\o{\overline}

\def\i0i{\int_0^\infty}

\def\op{\overline{p}}

\numberwithin{equation}{section}

\title{An improvement for the sharp Adams inequalities in bounded domains and whole space $\R^n$}
\author{Van Hoang Nguyen\footnote{
Institut de Math\'ematiques de Toulouse, Universit\'e Paul Sabatier, 118 Route de Narbonne, 31062 Toulouse c\'edex 09, France.}
}

\begin{document}
\maketitle

%\begin{center}
%\emph{(To my wife and my daughter)\text}
%\end{center}

%% Classification and key words; note that the 2010 classification is used:

\renewcommand{\thefootnote}{}

\footnote{Email: \href{mailto: Van Hoang Nguyen <van-hoang.nguyen@math.univ-toulouse.fr>}{van-hoang.nguyen@math.univ-toulouse.fr}}

\footnote{2010 \emph{Mathematics Subject Classification\text}: 26D10; 35B33; 46E35; 46E30.}

\footnote{\emph{Key words and phrases\text}: Moser--Trudinger inequality, Adams inequality, Lions Concentration--Compactness principle, Decreasing rearrangement, Best constant, Sobolev spaces.}

\renewcommand{\thefootnote}{\arabic{footnote}}
\setcounter{footnote}{0}

\begin{abstract}
We prove an improvement for the sharp Adams inequality in $W^{m,\frac nm}_0(\Om)$ where $\Om$ is a bounded domain in $\R^n$ inspired by Lions Concentration--Compactness principle for the sharp Moser--Trudinger inequality. Our method gives an alternative approach to a Concentration--Compactness principle in $W^{m,\frac nm}_0(\Om)$ recently established by do \'O and Macedo. Moreover, when $m$ is odd, we obtain an improvement for their result by finding the best exponent in this principle. Our approach also is successfully applied to whole space $\R^n$ to establish an improvement for the sharp Adams inequalities in $W^{m,\frac nm}(\R^n)$ due to Ruf, Sani, Lam, Lu, Fontana and Morpurgo. This type of improvement is still unknown, in general, except the special case $m=1$ due to do \'O, de Souza, de Medeiros and Severo. Our method is a further development for the method of $\check{\rm C}$erny, Cianchi and Hencl combining with some estimates for the decreasing rearrangement of a function in terms of the one of its higher order derivatives.

%The Adams--Moser--Trudinger inequality is the borderline case of the Sobolev inequality. A Concentration--Compactness principle of Lions's type for this inequality was recently established by do \'O and Macedo \cite{doOMacedo2014}. In this paper, we make a  \cite{Cerny2013} to give a new proof for this result of do \'O and Macedo. More details, our result improves their result when the order of gradient is odd, and also gives the best exponent in this Concentration--Compactness principle. Finally, we establish a Concentration--Compactness principle with the best exponent for the sharp Adams--Moser--Trudinger inequality in unbounded domains which was recently proved by Fontana and Morpurgo \cite{FM2015}.

%We give a new proof for the generalization of the Concentration--Compactness principle for the sharp Adams--Moser--Trudinger inequality in the Sobolev space $W_0^{m,\frac nm}(\Om)$ with $n\geq 2$, $m < n$ and $\Om$ is a bounded domain in $\R^n$ which was recently established by do \'O and Macedo \cite{doOMacedo2014}. Our result covers the one of do \'O and Macedo when $m$ is even, and improves their result when $m$ is odd. Our approach is a further development of the ideas 
\end{abstract}

\section{Introduction}
Let $\Om$ be a bounded domain (i.e., open subset) in $\R^n$, $n\geq 2$. The Sobolev embedding theorems assert that $W_0^{k,p}(\Om) \hookrightarrow L^q(\Om)$ where $k$ is a positive integer, $p \in (1, n/k),$ and $1 \leq q \leq np/(n-kp)$. Such theorems play an important role and are the central tools in many areas such as analysis, differential geometry, partial differential equations, calculus of variations, etc. However, in the limit case, $kp =n$, it is well-known that $W_0^{k,\frac nk}(\Om)\not\subset L^\infty(\Om)$. In this case, the Moser-Trudinger and Adams inequalities are the perfect replacement. The Moser-Trudinger inequality was established independently by ${\rm Yudovi\check{c}}$ \cite{Y1961}, ${\rm Poho\check{z}aev}$ \cite{P1965} and Trudinger \cite{T1967} which asserts the existence of a constant $\al > 0$ such that $W_0^{1,n}(\Om) \hookrightarrow L_\phi(\Om)$, where $L_\phi(\Om)$ is the Orlicz space determined by the Young function $\phi(t) = e^{\al |t|^{n/(n-1)}} -1$. Later, Moser sharpened this result by finding the best constant $\al$ in the embedding above. More precisely, he proved that for any $\al \leq \al_n := n^{\frac n{n-1}} \omega_{n}^{\frac1{n-1}}$ where $\omega_{n}$ denotes the volume of the unit ball in $\R^n$, there exists a constant $c_0 > 0$ such that

%\begin{Theorema} {\bf (Moser \cite{M1970})}
 %Let $\alpha_n = n^{\frac n{n-1}} \omega_{n}^{\frac1{n-1}}$ where $\omega_{n}$ denotes the volume of the unit ball in $\R^n$. There exists a constant $c_0 > 0$ such that
\begin{equation}\label{eq:Moser-Trudinger}
\sup\limits_{u \in C_0^\infty(\Om), \int_\Om |\nabla u|^n dx \leq 1}\frac1{|\Omega|} \int_\Om \exp\lt(\alpha |u|^{\frac n{n-1}}\rt) dx \leq c_0
\end{equation}
for any $\alpha\leq \alpha_n,$ and any bounded domain $\Om$ in $\R^n$. Moreover, the constant $\alpha_n$ is sharp in the sense that if $\alpha > \alpha_n$, then the supremum above will become infinity.
%\end{Theorema}

Moser--Trudinger inequality \eqref{eq:Moser-Trudinger} has played the important roles and has been widely applied in geometric analysis and PDEs, see for examples \cite{CY2003,LamLu2012a,LamLu2012b,LamLu2014,Shaw1987,TZ2000} and references therein. In recent years, it has been generalized in many directions, for instance, the singular Moser--Trudinger inequality \cite{AS2007}, or the sharp Moser--Trudinger inequality on domains of finite measure on the Heisenberg groups \cite{CohnLu2001,LamLuTang2012}, on spheres \cite{Bec1993}, on CR spheres \cite{CohnLu2001,CohnLu2004}, on the compact Riemannian manifolds \cite{Fontana1993,Li2005}, and on the hyperbolic spaces \cite{LT2013}. It was also extended to higher order of derivatives by Adams \cite{Adams1988}. To state Adams inequality, let us introduce some notations. We will use the symbol $\nabla^m$ with $m$ is a positive integer to denote the $m-$th order gradient of functions $u \in C^m(\Om)$, i.e.
\[
\nabla^m u =
\begin{cases}
\Delta^{\frac m2} u&\mbox{if $m$ is even}\\
\nabla\Delta^{\frac{m-1}2} u&\mbox{if $m$ is odd.}
\end{cases}
\]
We also use $\|\nabla^m u\|_p$ to denote the $L^p-$norm, $1\leq p < \infty$, of $\nabla^m u$. For $m < n$ and $1 \leq p < \infty$, we define the Sobolev space $W_0^{m,p}(\Om)$ as the completion of $C_0^\infty(\Om)$ under the norm $\|\nabla ^m u\|_p$ with $u\in C_0^\infty(\Om)$. In \cite{Adams1988}, Adams generalized the Moser--Trudinger inequality \eqref{eq:Moser-Trudinger} to the Sobolev spaces of higher order $W_0^{m,\frac nm}(\Om)$ as follows:
%\begin{Theoremb}{\bf (Adams \cite{Adams1988})}
Let $m$ be a positive integer less than $n$, there exists a constant $C(n,m) >0$ such that 
\begin{equation}\label{eq:Adams}
\sup\limits_{u\in W^{m,\frac nm}_0(\Om), \|\nabla^m u\|_{\frac nm} \leq 1}\frac{1}{|\Om|}\int_\Om \exp(\be(n,m) |u|^{\frac n{n-m}}) dx \leq C(n,m)
\end{equation}
where
\[
\beta(n,m) = 
\begin{cases}
\frac1{\om_n} \lt(\frac{\pi^{\frac n2} 2^m \Gamma\lt(\frac{m+1}2\rt)}{\Gamma\lt(\frac{n-m+1}2\rt)}\rt)^{\frac n{n-m}}&\mbox{if $m$ is odd}\\
\frac1{\om_n} \lt(\frac{\pi^{\frac n2} 2^m \Gamma\lt(\frac{m}2\rt)}{\Gamma\lt(\frac{n-m}2\rt)}\rt)^{\frac n{n-m}}&\mbox{if $m$ is even.}
\end{cases}
\]
Note that $\al_n =\beta(n,1)$. Moreover, the constant $\be(n,m)$ in \eqref{eq:Adams} is sharp in the sense that if we replace it by any $\beta > \be(n,m)$ then the supremum above will become infinite.
%\end{Theoremb}

The Adams approach to \eqref{eq:Adams} consists of some main steps. First, he represented $u$ in terms of $\nabla^m u$ via convolution with the Riesz potential, then applying O'Neil lemma \cite{O'Neil1963} to obtain a bound for the rearrangement function of $u$ via a the rearrangement function of $\nabla^m u$, and finally he applied Adams--Garsia lemma \cite{Adams1988} to obtain \eqref{eq:Adams}. Adams approach recently was used with some modifications to obtain the sharp Adams inequality in measure space by Fontana and Morpurgo \cite{FM2011}.

%The existence of extremal functions for Moser--Trudinger inequality \eqref{eq:Moser-Trudinger} was first proved by Carleson and Chang \cite{CC1986} when $\Om$ is the unit ball. They showed that the supremum in \eqref{eq:Moser-Trudinger} is achieved when $\alpha \leq \alpha_n$. This fact is very interesting and strange since the extremal functions in Sobolev inequalities do not exist for the bounded domains. Later, Struwe \cite{Struwe} obtained a sufficient condition for the existence of the extremal functions in the non-symmetric domains in $\R^2$. In \cite{Flucher1992} Flucher introduce another method, the conformal rearrangement, and derived an isoperimetric inequality which implies the existence of extremal functions to any smooth bounded domain in $\R^2$. At last, Lin \cite{Lin1996} generalized the existence of extremal function to any smooth bounded domain in $\R^n$, $n\geq 2$.

It was observed by Lions \cite{Lions1985} that the embedding $W^{1,n}_0(\Om) \hookrightarrow L_\phi(\Om)$ is not compact. He also proved in \cite{Lions1985} that except for "small weak neighborhoods of zero function" this embedding is compact by improving the best constant $\alpha_n$. His result now is known as a Concentration--Compactness principle. To state this principle, let us denote $\mathcal{M}(\overline{\Om})$ the space of all Radon measures on $\R^n$ whose support is $\overline{\Om}$, and $u^\sharp$ the spherically symmetric decreasing rearrangement of $u$ (see Section \S2 below for its definition). The Lions Concentration--Compactness principle states that

\begin{Theorema}{\bf (Lions \cite{Lions1985})} Let $\{u_j\}_j \subset W^{1,n}_0(\Om)$ such that $\|\nabla u_j\|_n \leq 1$, $u_j \rightharpoonup u$ in $W^{1,n}_0(\Om)$ and $|\nabla u_j|^n \rightharpoonup \mu$ in $\mathcal{M}(\overline{\Om})$. Then
\begin{description}
\item (i) If $u\equiv 0$ and $\mu =\de_{x_0}$ the Dirac measure concentrated at some point $x_0 \in \o{\Om}$, then up to a subsequence $e^{\alpha_n |u_j|^{\frac n{n-1}}} \rightharpoonup c\de_{x_0} + \mathcal{L}_n$ in $\mathcal{M}(\o{\Om})$ for some $c \geq 0$, and $\mathcal{L}_n$ denotes the Lebesgue measure on $\o{\Om}$.

\item (ii) If $u\equiv 0$ and $\mu$ is not Dirac measure, then there exist constant $p > 1$ and $C > 0$ depending only on $p$ and $\Om$ such that
\[
\sup\limits_{j\geq 1} \int_\Om e^{\al_n p |u_j|^{\frac n{n-1}}} dx \leq C.
\]

\item (iii) If $u\not \equiv 0$, then for any $p\in [1,\eta)$ with $\eta =(1 -\|\nabla u^\sharp\|_n^n)^{-\frac 1{n-1}}$, there exists constant $C$ depending only on $p$ and $\Om$ such that
\[
\sup\limits_{j\geq 1} \int_\Om e^{\al_n p |u_j|^{\frac n{n-1}}} dx \leq C.
\]
\end{description}
\end{Theorema}
It is clear that Lions Theorem gives more informations than Moser--Trudinger inequality when $u_j \rightharpoonup  u$ weakly in $W^{1,n}_0(\Om)$ with $u\not\equiv 0$. We remark that the upper bound $\eta$ for the value of $p$ in the case $(iii)$ of Lions Theorem is not sharp. In \cite{Cerny2013}, $\check{\rm C}$erny, Cianchi and Hencl presented a new approach to Lions Theorem, and yield a sharp upper bound for these values of $p$. Their resut reads as follows.
\begin{Theoremb}{\bf (${\bf \check{{\rm {\bf C}}}}$erny, Cianchi and Hencl \cite{Cerny2013})}
Under the same assumptions as in the case $(iii)$ of Lions Theorem, define $P =(1 -\|\nabla u\|_n^n)^{-\frac 1{n-1}}$, then for any $p \in [1,P)$ there exists constant $C$ depending only on $p$ and $\Om$ such that
\[
\sup\limits_{j\geq 1} \int_\Om e^{\al_n p |u_j|^{\frac n{n-1}}} dx \leq C.
\]
Moreover, the upper bound $P$ for $p$ is sharp.
\end{Theoremb}
Note that one always has $\eta \leq P$ in general because of P\'olya--Szeg\"o principle, and the inequality is strict unless $u$ has a very special form (see \cite{BZ}).

It is also true that the embedding $W^{m,\frac nm}_0(\Om) \hookrightarrow L_\phi(\Om)$ with $\phi(t) = e^{\beta(n,m) |t|^{\frac n{n-m}}}-1$ is not compact. Motivated by Lions Theorem, do \'O and Macedo studied the compactness of this embedding in \cite{doOMacedo2014} and established a Concentration-Compactness principle of Lions type for the sharp Adams inequality \eqref{eq:Adams} by improving the best constant $\beta(n,m)$ in this inequality. For a function $u\in W^{m,\frac nm}_0(\Om)$, $\|\nabla^m u\|_{\frac nm} \leq 1$ let us denote
\[
\eta_{n,m}(u) =
\begin{cases}
(1 -\|\nabla^m u\|_{\frac nm}^{\frac nm})^{-\frac m{n-m}}&\mbox{if $m$ is even}\\
(1 -\|\nabla(\Delta^{\frac{m-1}2}u)^\sharp\|_{\frac nm}^{\frac nm})^{-\frac{m}{n-m}}&\mbox{if $m$ is odd.}
\end{cases}
\]
\begin{Theoremc}{\bf (do \'O and Macedo \cite{doOMacedo2014})}
Let $m$ be a positive integer with $m < n$ and $n/m \geq 2n / (n+2)$. Let $\{u_j\}_j \subset W^{m,\frac nm}_0(\Om)$ such that $\|\nabla^m u_j\|_{\frac nm} \leq 1$, and $u_j\rightharpoonup u$ weakly in $W^{m,\frac nm}_0(\Om)$ and $|\nabla^m u_j|^{\frac nm} \rightharpoonup \mu$ in $\mathcal{M}(\o{\Om})$. Then,
\begin{description}
\item (i) If $u\equiv 0$ and $\mu =\de_{x_0}$ the Dirac measure concentrated at some point $x_0 \in \o{\Om}$, then up to a subsequence $e^{\beta(n,m) |u_j|^{\frac n{n-m}}} \rightharpoonup c\de_{x_0} + \mathcal{L}_n$ in $\mathcal{M}(\o{\Om})$ for some $c \geq 0$.

\item (ii) If $u\equiv 0$ and $\mu$ is not Dirac measure, then there exist constant $p > 1$ and $C > 0$ depending only on $p$ and $\Om$ such that
\[
\sup\limits_{j\geq 1} \int_\Om e^{\beta(n,m) p |u_j|^{\frac n{n-m}}} dx \leq C.
\]

\item (iii) If $u\not \equiv 0$, then for any $p\in [1,\eta_{n,m}(u))$, there exists constant $C$ depending only on $p$ and $\Om$ such that
\begin{equation}\label{eq:doOMacedo}
\sup\limits_{j\geq 1} \int_\Om e^{\beta(n,m) p |u_j|^{\frac n{n-m}}} dx \leq C.
\end{equation}
\end{description}
\end{Theoremc}
It is worth to mention here that in the case $n =2m$, one can prove \eqref{eq:doOMacedo} for any $1 \leq p < (1 -\|\nabla^m u\|_2^2)^{-1}$ by using the same argument in \cite{Lions1985} exploiting the Hilbert structure of the space $W^{m,2}_0(\Om)$. This fact and the result of $\check{\rm C}$erny, Cianchi and Hencl suggest us that in the case $m$ odd, the inequality \eqref{eq:doOMacedo} also holds for any $1 \leq p < (1-\|\nabla^m u\|_{\frac nm}^{\frac nm})^{-\frac m{n-m}}$. We remark that the proof of do \'O and Macedo given in \cite{doOMacedo2014} follows the method of Lions \cite{Lions1985} based on the symmetrization argument. More precisely, they use the symmetrization argument over function $\Delta^{\frac{m-1}2}u$ when $m$ is odd, the Talenti comparison principle \cite{Tal1976} and maximum principle to prove \eqref{eq:doOMacedo}. This approach prevents us to reach the value $(1-\|\nabla^m u\|_{\frac nm}^{\frac nm})^{-\frac m{n-m}}$ of $\eta_{n,m}(u)$ when $m$ is odd. In this paper, by making a further development for the method of $\check{\rm C}$erny, Cianchi and Hencl, we give another proof of \eqref{eq:doOMacedo}. We also prove that in the case $m$ odd, we can attain the value $(1-\|\nabla^m u\|_{\frac nm}^{\frac nm})^{-\frac m{n-m}}$ of $\eta_{n,m}(u)$. This improves the previous result of do \'O and Macedo when $m$ is odd because of the P\'olya--Szeg\"o principle. For $u\in W^{m,\frac nm}_0(\Om)$, $\|\nabla^m u\|_{\frac nm} \leq 1$, let us denote
\[
P_{n,m}(u) = 
\begin{cases}
(1-\|\nabla^m u\|_{\frac nm}^{\frac nm} dx)^{-\frac{m}{n-m}}&\mbox{if $\|\nabla^m u\|_{\frac nm} < 1$}\\
\infty &\mbox{if $\|\nabla^m u\|_{\frac nm} =1$}.
\end{cases}
\]
We will prove the following result.
\begin{theorem}\label{Maintheorem}
Assume that $m$ is an integer less than $n$. Under the same assumptions in the case $(iii)$ of do \'O and Macedo Theorem, we have
\begin{equation}\label{eq:CCNgoNguyen}
\sup_{j} \int_{\Om} e^{p \beta(n,m) |u_j|^{\frac n{n-m}}} dx < \infty,
\end{equation}
for any $p\in [1, P_{n,m}(u))$. Moreover, the upper bound $P_{n,m}(u)$ is sharp.
\end{theorem}
%Because of the P\'olya--Szeg\"o principle, we have $P_{n,m}(u) \geq \eta_{n,m}(u)$ when $m$ is odd, hence our result \eqref{eq:CCNgoNguyen} is better than the one of do \'O and Macedo \eqref{eq:CCdoOMacedo} in the case $m$ odd. Our method to prove Theorem \ref{Maintheorem} goes back ${\rm \check{C}}$erny, Cianchi and Hencl \cite{Cerny2013} where they improved the concentration-compactness principle due to Lions \cite{Lions1985}. More precisely, they proved Theorem \ref{Maintheorem} for $m=1$. 

The Moser--Trudinger inequality was first extended to unbounded domains by Cao \cite{Cao} in $\R^2$, and then for any dimension $n\geq 2$ by do \'O \cite{doO}. Later, Adachi and Tanaka \cite{AT1999} established a sharp version for the result of Cao and do \'O. Note that all these inequalities were proved by taking the supremum with respect to the Dirichlet norm of gradient and assumed in some sense a subcritical growth $e^{\alpha |u|^{\frac n{n-1}}}$ with $\alpha < \al_n$. More precisely, we have the following inequality \cite{AT1999}: for any $\al < \alpha_n$, there exists a constant $C(n,\al)$ depending only on $n$ and $\al$ such that
\begin{equation}\label{eq:AdachiTanaka}
\sup\limits_{u\in W^{1,n}(\R^n), \|\nabla u\|_n \leq 1} \int_{\R^n} \Phi_{n,1}(\al |u|^{\frac n{n-1}}) dx \leq C(n,\al),
\end{equation}
where $\Phi_{n,1}(t) =e^t - \sum_{k=0}^{n-2} \frac {t^k}{k!}$. Moreover, the constant $\alpha_n$ is sharp in the sense that the supremum above will become infinite if $\al \geq \alpha_n$. In recent paper \cite{LamLuZhang2015}, Lam, Lu and Zhang proved that the constant $C(n,\al)$ satisfies $C(n,\al) \leq C(n)/(\al_n -\al)$ for some constant $C(n)$ depending only on $n$. Recently, Li and Ruf \cite{LiRuf2008,Ruf2005} showed that the Moser--Trudinger inequality can be extended to any bounded domains (and thus to all of $\R^n$) with the critical exponent $\alpha_n$ if we replace the Dirichlet norm of gradient by the full Sobolev norm $\|u\|_{n,1} = (\|u\|_n^n + \|\na u\|_n^n)^{\frac 1n}$
\begin{equation}\label{eq:LiRuf}
\sup\limits_{u\in W^{1,n}(\R^n), \|\nabla u\|_{n,1} \leq 1} \int_{\R^n} \Phi_{n,1}(\al_n |u|^{\frac n{n-1}}) dx  <\infty.
\end{equation}
It is very strange and interesting that the inequalities \eqref{eq:AdachiTanaka} and \eqref{eq:LiRuf} are equivalent, in general, as shown in \cite{LamLuZhang2015}.

Inspired by the Concentration--Compactness principle due to Lions (Theorem {\bf A} above), do \'O, de Souza, de Medeiros and Severo \cite{doO2014} recently proved the following result.
\begin{Theoremd}{\bf (do \'O, de Souza, de Medeiros, and Severo \cite{doO2014})}
Let $\{u_j\}_j \subset W^{1,n}(\R^n)$ such that $\|u_j\|_{n,1} =1$ and $u_j \rightharpoonup u\not\equiv 0$ in $W^{1,n}(\R^n)$, denote $Q_{n,1}(u) = (1 -\|u\|_{n,1}^n)^{-1/(n-1)}$, then 
\[
\sup\limits_{j\in \N} \int_{\R^n} \Phi_{n,1}(p \al_n |u_j|^{\frac n{n-1}}) dx < \infty,
\]
for any $p \in [1,Q_{n,1}(u))$. Moreover, the upper bound $Q_{n,1}(u)$ for $p$ is sharp.
\end{Theoremd}
It is clear that Theorem {\bf D} improves the Moser--Trudinger inequality \eqref{eq:LiRuf}. It completes some results in \cite{LamLu2012a,Yang2012a} where the authors established a singular versions of Theorem {\bf D} under the additional assumptions that $\nabla u_j$ converges a.e to $\nabla u$ in $\R^n$. Obviously, this assumption is restrictive and we emphasize that it was crucial for the arguments in \cite{LamLu2012a,Yang2012a} which follows the lines of Lions \cite{Lions1985}. The proof of Theorem {\bf D} given in \cite{doO2014} follows the lines of $\check{\rm C}$erny, Cianchi and Hencl \cite{Cerny2013}. In the dimension two, the proof of Theorem {\bf D} is easy by exploiting the Hilbert structure of the space $W^{1,2}(\R^2)$ (see \cite{doO2008}).

Extending the Adams inequality \eqref{eq:Adams} to unbounded domains in $\R^n$ is an interesting problem. This problem was first done by Ruf and Sani \cite{RufSani2013} for $m$ even, and then by Lam and Lu \cite{LamLu2012d} for an arbitrary $m$. To state their inequality, let us denote a norm $\| \cdot\|$, for $u \in W^{m,\frac nm}(\R^n)$, by
\[
\|u\|_{m,\frac nm} =
\begin{cases}
\|(-\Delta + I)^{\frac m2} u\|_{\frac nm} &\mbox{if $m$ is even}\\
\lt(\|(-\Delta + I)^{\frac{m-1}2} u\|_{\frac nm}^{\frac nm} + \|\nabla (-\Delta + I)^{\frac{m-1}2} u\|_{\frac nm}^{\frac nm}\rt)^{\frac mn} &\mbox{if $m$ is odd.}
\end{cases}
\]
Then the following inequality was proved in \cite{RufSani2013,LamLu2012d}, 
\begin{equation}\label{eq:RufSaniLamLu}
\sup\limits_{u \in W^{m,\frac nm}(\R^n),\, \|u\|_{m,\frac nm} \leq 1}\int_{\R^n} \Phi_{n,m}(\beta(n,m) |u|^{\frac n{n-m}}) dx < \infty,
\end{equation}
where $\Phi_{n,m}(t) =e^t -\sum_{j=0}^{j_{\frac nm}- 2} \frac{t^j}{j!}$ with
\[
j_{\frac nm} = \min\lt\{j \in \N \, :\, j \geq \frac nm\rt\} \geq \frac nm.
\]
The constant $\beta(n,m)$ in \eqref{eq:RufSaniLamLu} is sharp in the sense that the supremum in \eqref{eq:RufSaniLamLu} will become infinite if we replace $\beta(n,m)$ by any larger constant. We refer the reader to the paper \cite{LamLu2013} for a sharp Adams type inequality of fractional order $\alpha \in (0,n)$ where a rearrangement--free argument was found. 

It was asked in \cite{doOMacedo2014} that proving a Lions type improvement of the Adams inequality \eqref{eq:RufSaniLamLu}, like Theorem {\bf D}, is an interesting question. We remark that this question is still open except a few cases in dimension four given in \cite{Yang2012} where the author exploited the Hilbert structure of the space $W^{2,2}(\R^4)$. Our next result provides such an improvement for the Adams inequality \eqref{eq:RufSaniLamLu} in full generality.

\begin{theorem}\label{Maintheoremfull}
Let $\{u\}_j \subset W^{m,\frac nm}(\R^n)$ such that $\|u_j\|_{m,\frac nm} \leq 1$, $u_j$ converges weakly to $u\not\equiv 0$ in $W^{m,\frac nm}(\R^n)$, then 
\begin{equation}\label{eq:CCunboundedfull}
\sup_{j} \int_{\R^n} \Phi_{n,m}(p \beta(n,m) |u_j|^{\frac n{n-m}}) dx < \infty,
\end{equation}
for any $p\in [1, Q_{n,m}(u))$, with
\[
Q_{n,m}(u) =
\begin{cases}
(1 -\|u\|_{m,\frac nm}^{\frac nm} )^{-\frac m{n-m}}&\mbox{if $\|u\|_{m,\frac nm} < 1$}\\
\infty &\mbox{if $\|u\|_{m,\frac nm} =1$.}
\end{cases}
\]
\end{theorem}

We next introduce the Ruf norm on $W^{m,\frac nm}(\R^n)$ by
\[
\|u\| = \lt(\|u\|_{\frac nm}^{\frac nm} + \|\nabla^m u\|_{\frac nm}^{\frac nm}\rt)^{\frac mn}, \quad u\in W^{m,\frac nm}(\R^n).
\]
Recently, Fontana and Morpurgo established in \cite{FM2015} a sharp Adams inequality under the \emph{Ruf condition} $\|u\| \leq 1$. Their inequality reads as follows
\begin{equation}\label{eq:FontanaMorpurgo}
\sup\limits_{u \in W^{m,\frac nm}(\R^n),\, \|u\| \leq 1}\int_{\R^n} \Phi_{n,m}(\beta(n,m) |u|^{\frac n{n-m}}) dx < \infty.
\end{equation}
The constant $\beta(n,m)$ in \eqref{eq:RufSaniLamLu} is sharp in the sense that the supremum in \eqref{eq:RufSaniLamLu} will become infinite if we replace $\beta(n,m)$ by any larger constant. We emphasize that the \emph{Ruf condition} is in some sense minimal, in regarding to the number of derivatives to obtain the sharp Adams inequality. Such a type of this inequality when $m=2$ was recently proved by Lam and Lu \cite{LamLu2013} by the domain decomposition method which is completely with the one of Fontana and Morpurgo.

%Our second aim of this paper is to establish the concentration-compactness principle for sharp Adams inequality in unbounded domains. 

%Let $\Om$ be an unbounded domain of $\R^n$, we denote $W^{m,\frac nm}(\R^n)$ the completion of $C_0^\infty(\Om)$ under the norm $\lt(\int_\Om |u|^{\frac nm} dx + \int_\Om |\nabla^m u|^{\frac nm} dx \rt)^{\frac mn},$ with $ u\in C_0^\infty(\Om)$. For $u \in W^{m,\frac nm}_0(\Om)$ such that  $\int_\Om |u|^{\frac nm} dx + \int_\Om |\nabla^m u|^{\frac nm} dx \leq 1$, we denote

Our next result gives us an improvement of the sharp Adams inequality \eqref{eq:FontanaMorpurgo} under the \emph{Ruf condition} by improving the constant $\beta(n,m)$ in the light of the Lions Concentration--Compactness principle. More precisely, we prove the following theorem.
\begin{theorem}\label{Maintheorem2}
Let $\{u\}_j\subset W^{m,\frac nm}(\R^n)$ such that $\|u_j\| \leq 1$, $u_j$ converges weakly to $u\not\equiv 0$ in $W^{m,\frac nm}(\R^n)$, then 
\begin{equation}\label{eq:CCunbounded}
\sup_{j} \int_{\R^n} \Phi_{n,m}(p \beta(n,m) |u_j|^{\frac n{n-m}}) dx < \infty,
\end{equation}
for any $p\in [1, R_{n,m}(u))$ with 
\[
R_{n,m}(u) =
\begin{cases}
(1 -\|u\|^{\frac nm})^{-\frac m{n-m}}&\mbox{if $\|u\| < 1$}\\
\infty &\mbox{if $\|u\| =1$.}
\end{cases}
\]
\end{theorem}
Let us make some comments on the proofs of our main Theorems \ref{Maintheorem}, \ref{Maintheoremfull} and \ref{Maintheorem2}. Our approach is based on the method of $\check{\rm C}$erny, Cianchi and Hencl in \cite{Cerny2013} where the authors improved the Lions Concentration--Compactness principle for the sharp Moser--Trudinger inequality by finding the best exponent in this principle. To apply the method of $\check{\rm C}$erny, Cianchi and Hencl, we need to establish some estimates for the decreasing rearrangement of a function in terms of the one of its derivatives in higher order. These estimates generalize a recent estimate of Masmoudi and Sani \cite{MS2014} to higher order derivatives, and seem to be new. It is worth to mention here that our method used in this paper recently was applied in \cite{NgoNguyen2016} to establish a Concentration--Compactness principle of Lions type in the Hyperbolic space. This result extends the recent results of Karmakar \cite{Kar2015} for $m=1,2$ to higher order of derivatives (i.e, $m \geq 3$), and gives a alternative proof for the result of Karmakar when $m=2$.

We finish this introduction by recalling an elementary inequality which is used frequently in this paper. Given $p >1$ and $\ep >0$ then the following estimate holds
\begin{equation}\label{eq:elementary}
(a +b)^p \leq (1+\ep) a^p + C_\ep b^p,\quad\forall\, a,\, b > 0
\end{equation}
with $C_\ep = (1 - (1+\ep)^{-\frac1{p-1}})^{1-p}$.

The rest of this paper is organized as follows. In the next section \S2, we recall the notions of the decreasing rearrangement and the spherically decreasing rearrangement functions of a given function. We also prove in this section some useful estimates involving the rearrangement of solutions of the polyharmonic equations which will be crucial in the proofs of our main Theorems. Section \S3 is devoted to prove Theorem \ref{Maintheorem} and the sharpness of this theorem. The proof of Theorem \ref{Maintheoremfull} and Theorem \ref{Maintheorem2} are given in Section \S4. %In the section \S5, we show how a minor modification of the proof of Theorem \ref{Maintheoremfull} proves Theorem \ref{Maintheorem2}. We finish this paper by proving Corollary \ref{outsidezerofunction} in the Section \S6.

\section{Some preliminaries}
\subsection{Rearrangement}
We start this section by recalling the notion of the decreasing rearrangement function of a given function $u$ defined in a subset of $\R^n$. Let $\Omega\subset\R^n$ be a measurable set, we denote $\Omega^\sharp$ the open ball $B_R$ centered at origin of radius $R>0$ such that 
\[
|\Omega|= |B_R| = R^n \omega_n.
\]
Let $u: \Omega\to \R$ be a real-valued measurable function in $\Omega$. Then the distribution function of $u$ is the function $\mu_u:[0,\infty)\to [0,\infty]$ defined as 
\[
\mu_u(s) = |\{x\in \Omega \, :\, |u(x)| > s\}|,\quad s \geq 0,
\]
and the decreasing rearrangement of $u$ is the right-continuous, non-increasing function $u^*:[0,\infty)\to [0,\infty]$ which is equimeasurable with $u$, namely
\[
u^*(t) = \sup\{s \geq 0\, :\, \mu_u(s) >t\},\quad t \geq 0.
\]
Note that the support of  $u^*$ satisfies ${\rm supp}\, (u^*)\subset [0,|\Omega|]$. Since $u^*$is non-increasing, the maximal function $u^{**}$ of the rearrangement $u^*$, defined by
\[
u^{**}(t) = \frac1t \int_0^t u^*(s) ds,\quad t\geq 0,
\]
is also non--increasing and $u^* \leq u^{**}$. Moreover, we have the following inequality.
\begin{proposition}\label{normrearrangement}
If $u\in L^p(\Omega)$ with $1 < p< \infty$ and $\frac1p + \frac1{p'} =1$, then
\[
\left(\int_0^\infty (u^{**}(t))^p dt\right)^{\frac1p} \leq p' \left(\int_0^\infty (u^*(t))^pdt\right)^{\frac1p}.
\]
In particular, if ${\rm supp}\, (u )\subset \Omega$ with $\Omega$ a domain in $\R^n$, then
\[
\left(\int_0^{|\Omega|}(u^{**}(t))^p dt\right)^{\frac1p} \leq p' \left(\int_0^{|\Omega|} (u^*(t))^pdt\right)^{\frac1p}.
\]
\end{proposition}

Finally, we will denote by $u^\sharp: \Omega^\sharp \to [0,\infty]$ the spherically symmetric decreasing rearrangement of $u$
\[
u^\sharp(x) = u^*(\omega_n |x|^n),\quad x\in \Omega^\sharp.
\]
The function $u^\sharp$ has the same distribution of $u$, hence for any $p\in [1,\infty)$, we have
\[
\int_{\Om^\sharp} |u^\sharp(x)|^p dx = \int_{\Om} |u(x)|^p dx.
\]
Moreover, if $u\in W^{1,p}_0(\Om)$ then $u^\sharp \in W^{1,p}_0(\Om^\sharp)$ and
\[
\int_{\Om^\sharp} |\nabla u^\sharp|^p dx \leq \int_{\Om} |\nabla u|^p dx,
\]
by P\'olya--Szeg\"o principle.

\subsection{Some useful inequalities involving the rearrangement}
Let us start this section by recalling an useful and interesting result of Masmoudi and Sani \cite{MS2014}. Let $f\in L^2(\Omega)$, we consider the following Dirichlet problem
\begin{equation}\label{eq:Dirichlet}
\begin{cases}
-\Delta u = f&\mbox{in $\Omega$}\\
u = 0&\mbox{on $\partial \Omega$}.
\end{cases}
\end{equation}
If $u \in W^{1,2}_0(\Omega)$ be the unique weak solution to \eqref{eq:Dirichlet}, then 
\begin{equation}\label{eq:MSrearrangement}
u^*(t_1) -u^*(t_2) \leq \frac{1}{n^2 \omega_n^{\frac2n}}\int_{t_1}^{t_2} \frac{f^{**}(t)}{t^{1-\frac2n}} dt
\end{equation}
for any $0 < t_1 < t_2 < |\Omega|$. Since $u^*(|\Om|) =0$, then by integration by parts, we obtain that
\begin{equation}\label{eq:uhaisao}
u^{**}(t) \leq \frac 1{(n\om_n^{1/n})^2} \lt(\int_t^{|\Omega|} \frac{f^{**}(s)}{s^{1 -\frac 2n}} ds + \frac1t\int_0^t f^{**}(s) s^{\frac2n} ds\rt),
\end{equation}
for any $0 < t < |\Om|$, here we use the simple fact 
\[
\lim_{t\to 0^+} t \int_t^{|\Omega|} \frac{f^{**}(s)}{s^{1 -\frac 2n}} ds = 0.
\]
The inequality \eqref{eq:MSrearrangement} is a crucial ingredient in the proof of the Adams inequality with exact growth condition in \cite{MS2014} for dimension $4$ and in \cite{LTZ2015} for any dimension $n\geq 3$. Our first aim of this section is to generalize the inequality \eqref{eq:MSrearrangement} to higher order of derivatives. Such a generalization is as follows

\begin{proposition}\label{keypropo}
Let $u\in C_0^\infty(\R^n)$ and let $k$ be a positive integer less than $\frac n2$, denote $f =(-\Delta)^k u \in C_0^\infty(\R^n)$,  then there exists a constant $C(n,k) > 0$ depending only on $n,k$ such that 
\begin{equation}\label{eq:keyinproof}
u^*(t_1) - u^*(t_2) \leq \frac{c(n,k)}{(n\omega_n^{\frac1n})^{2k}} \int_{t_1}^{t_2} \frac{f^{**}(s)}{s^{1 -\frac{2k}n}} ds + C(n,k)\|f\|_{\frac n{2k}},
\end{equation}
where 
\begin{equation*}
c(n,k) = 
\begin{cases}
1&\mbox{if $k=1$}\\
\frac{n^{2(k-1)}}{2^{k-1}(k-1)! \Pi_{j=1}^{k-1}(n-2j)}&\mbox{if $k\geq 2$}
\end{cases}
\end{equation*}
for any $0< t_1 < t_2 < \infty$.
\end{proposition}
We make a remark here that the estimate \eqref{eq:keyinproof} seems to be not good because of the appearance of the term $C(n,k) \|f\|_{\frac n{2k}}$ in its right hand side. However, this estimate is enough for us to prove Theorems \ref{Maintheorem} and \ref{Maintheorem2}. It is worth to emphasize that if this term is dropped in \eqref{eq:keyinproof}, we can obtain the Adams inequality for higher derivatives $k>1$ with exact growth condition under the condition $\|\Delta^k u\|_{\frac n{2k}} \leq 1$. Let us go to the proof of Proposition \ref{keypropo}, for simplicity, we define the function $g$ in $\R_+\times \R_+$ by
\[
g(t,s) = 
\begin{cases}
s^{-1 + \frac2n}&\mbox{if $ t\leq s$}\\
t^{-1}s^{\frac2n}&\mbox{if $t >s$.}
\end{cases}
\]
We next define consecutively the sequence of functions $g_i(t,s)$ for $i =1,2,\ldots$ by
\[
g_1:=g,\quad g_i(t,s) = \int_0^\infty g_{i-1}(t,s_1) g(s_1,s) ds_1,\quad i\geq 2.
\]
Then we have the following estimates for $g_i$.
\begin{lemma}\label{boundsforg}
There exists a constant $C >0$ depending only on $n$ and $k$ such that
\begin{equation*}
g_i(t,s) \leq 
\begin{cases}
C s^{-1 + \frac{2i}n} &\mbox{if $t \leq s$}\\
C t^{-1 +\frac{2(i-1)}n} s^{\frac 2n}&\mbox{if $t > s$,}
\end{cases}
\end{equation*}
for any $i = 1,2,\ldots,k-1$.
\end{lemma}
\begin{proof}
We denote, in this proof, by $C$ a positive constant which depends only on $n$ and $k$, and can be changed from line to line. We prove by induction argument. The conclusion is obviously true for $i=1$. We next suppose that the conclusion is true for $1\leq i < k-1$, we then have to prove that it is true for $i+1$. Indeed, if $t < s$, an easy computation shows that
\begin{equation*}
g_{i}(t,s_1)g(s_1,s) \leq 
\begin{cases}
Ct^{-1+\frac{2(i-1)}n} s^{-1+\frac2n} s_1^{\frac2n}&\mbox{if $s_1 \leq t$}\\
Cs^{-1+\frac2n} s_1^{-1+\frac{2i}n}&\mbox{if $t< s_1 \leq s$}\\
Cs_1^{-2 + \frac{2i}n}s^{\frac2n}&\mbox{if $s_1 > s$,}
\end{cases}
\end{equation*}
which immediately implies that
\[
g_{i+1}(t,s) = \int_0^\infty g_i(t,s_1) g(s_1,s) ds_1 \leq C s^{-1+\frac{2(i+1)}n},
\]
when $t < s$.

If $t > s$, as above we have
\begin{equation*}
g_{i}(t,s_1)g(s_1,s) \leq 
\begin{cases}
Ct^{-1+\frac{2(i-1)}n} s^{-1+\frac2n} s_1^{\frac2n}&\mbox{if $s_1 \leq s$}\\
Ct^{-1+\frac{2(i-1)}n} s^{\frac2n} s_1^{-1+\frac2n}&\mbox{if $s< s_1 \leq t$}\\
Cs_1^{-2 + \frac{2i}n}s^{\frac2n}&\mbox{if $s_1 > t$,}
\end{cases}
\end{equation*}
which then easily implies 
\[
g_{i+1}(t,s) = \int_0^\infty g_i(t,s_1) g(s_1,s) ds_1 \leq C t^{-1+\frac{2i}n} s^{\frac2n},
\]
when $t > s$. This finishes the proof of Lemma \ref{boundsforg}.
\end{proof}
We also need the following result in the proof of Proposition \ref{keypropo}.
\begin{lemma}\label{derivativeofg}
For any $i=2,3,\ldots, k-1$, we have
\[
\partial_t g_i(t,s)=-\frac1{t^2} \int_0^t s_1^{\frac2n} g_{i-1}(s_1,s) ds_1.
\]
\end{lemma}
\begin{proof}
It is evident that
\begin{align*}
g_i(t,s) &= \int_0^\infty g(t,s_1) g_{i-1}(s_1,s) ds_1\\
&= \frac1t \int_0^t s_1^{\frac2n} g_{i-1}(s_1,s) ds_1 + \int_t^{\infty} s_1^{-1+\frac2n} g_{i-1}(s_1,s) ds_1,
\end{align*}
which then implies our desired result.
\end{proof}

\emph{Proof of Proposition \ref{keypropo}:} If $k=1$ then \eqref{eq:keyinproof} is an easy consequence of \eqref{eq:MSrearrangement}, hence we only prove \eqref{eq:keyinproof} for $k\geq 2$ in the rest of proof. Denote $u_i = (-\Delta)^i u$ for $i =0,1,\ldots,k$, then $u_k =f$. It follows from \eqref{eq:uhaisao} that
%\begin{equation}\label{eq:MSapply}
%-(u_i^*)'(t) \leq \frac1{n^2 \omega_n^{\frac2n}} \frac{u_{i+1}^{**}(t)}{t^{1 -\frac2n}},\quad 0< t< |\Omega|
%\end{equation}
%for $i=0,1,\ldots, k-1$. Then
\begin{align*}
u_{i}^{**}(t)&\leq \frac1{n^2 \omega_n^{\frac2n}}\left(\int_t^{\infty} \frac{u_{i+1}^{**}(s)}{s^{1 -\frac 2n}} ds + \frac1t\int_0^t u_{i+1}^{**}(s) s^{\frac2n} ds\right)\\
& \leq \frac1{n^2 \omega_n^{\frac2n}} \int_0^\infty g(t,s)u_{i+1}^{**}(s) ds,
\end{align*}
for any $0 < t< \infty$, and for $i =0,1,\ldots, k-1$. Using the latter inequality and Fubini theorem hence implies that
\[
(n\omega_n^{\frac1n})^{2(k-1)}u_1^{**}(t) \leq \int_0^\infty f^{**}(s) g_{k-1}(t,s) ds,\quad \forall\, t\in (0, |\Om|).
\]
Consequently, for any $0<t_1< t_2 < \infty$, we obtain
\begin{align}\label{eq:11111111}
(n\omega_n^{\frac1n})^{2k}(u^{*}(t_1)-u^*(t_2))& \leq \int_{t_1}^{t_2} \lt(\int_0^\infty f^{**}(s) g_{k-1}(t,s) ds\rt) t^{\frac2n-1}dt \notag\\
&= \int_0^\infty f^{**}(s)\left(\int_{t_1}^{t_2} g_{k-1}(t,s) t^{\frac2n-1} dt\right)ds
\end{align}
We next claim that there exist the functions $F_i(t_1,t_2,s)$ for any $i =1, 2,\ldots, k-1$ and $0< t_1 < t_2< \infty$  and constant $C >0$ depending only on $n,k$ such that 
\begin{equation}\label{eq:recurrentclaimforg}
\int_{t_1}^{t_2} g_{k-i}(t,s) t^{\frac{2i}n-1} dt \leq \frac{n^2}{2i(n-2i)} \int_{t_1}^{t_2} g_{k-i-1}(t,s) t^{\frac{2(i+1)}n-1} dt + F_i(t_1,t_2,s)
\end{equation}
for $i =1,2,\ldots, k-2$ if $k>2$, and
\begin{equation}\label{eq:forg1}
\int_{t_1}^{t_2} g(t,s) t^{\frac{2(k-1)}n-1} dt \leq \frac{n}{2(k-1)} \int_{t_1}^{t_2} s^{\frac2n} t^{\frac{2(k-1)}n-2} \chi_{\{s < t\}} dt + F_{k-1}(t_1,t_2,s)
\end{equation}
with
\begin{equation}\label{eq:boundnormclaim}
\left(\int_0^\infty F_i(t_1,t_2,s)^{\frac{n}{n-2k}} ds\right)^{\frac{n-2k}n} \leq C,
\end{equation}
Indeed, for $i=1,2,\ldots, k-2$ (if $k>2$), we have
\[
\int_{t_1}^{t_2} g_{k-i}(t,s) t^{\frac{2i}n-1} dt = \frac n{2i}\int_{t_1}^{t_2} g_{k-i}(t,s) dt^{\frac{2i}n}.
\]
Using integration by parts and Lemma \ref{derivativeofg}, we obtain 
\begin{align*}
\int_{t_1}^{t_2} g_{k-i}(t,s) dt^{\frac{2i}n} &= g_{k-i}(t_2,s) t_2^{\frac{2i}n} -g_{k-i}(t_1,s)t_1^{\frac{2i}n} - \int_{t_1}^{t_2} \partial_t g_{k-i}(t,s) t^{\frac{2i}n} dt\\
&= g_{k-i}(t_2,s) t_2^{\frac{2i}n} -g_{k-i}(t_1,s)t_1^{\frac{2i}n} + \int_{t_1}^{t_2} \frac1{t^2}\int_0^t s_1^{\frac2n} g_{k-i-1}(s_1,s) ds_1 t^{\frac{2i}n} dt\\
&= g_{k-i}(t_2,s) t_2^{\frac{2i}n}-g_{k-i}(t_1,s)t_1^{\frac{2i}n} -\frac{n}{n-2i} \int_{t_1}^{t_2}\int_0^t s_1^{\frac2n} g_{k-i-1}(s_1,s) ds_1 dt^{\frac{2i}n-1}\\
&=g_{k-i}(t_2,s) t_2^{\frac{2i}n}-g_{k-i}(t_1,s)t_1^{\frac{2i}n} -\frac n{n-2i}t_2^{\frac{2i-n}n}\int_0^{t_2}s_1^{\frac2n}g_{k-i-1}(s_1,s) ds_1\\
&\,+\frac{n}{n-2i} t_1^{\frac{2i-n}n}\int_0^{t_1}s_1^{\frac2n}g_{k-i-1}(s_1,s) ds_1 + \frac n{n-2i}\int_{t_1}^{t_2} g_{k-i-1}(t,s) t^{\frac{2(i+1)-n}n} dt \\
%&\leq g_{k-i}(t_2,s) t_2^{\frac{2i}n}+\frac n{n-2i}\int_{t_1}^{t_2} g_{k-i-1}(t,s) t^{\frac{2(i+1)-n}n} dt \\
%&\quad + \frac{n}{n-2i} t_1^{\frac{2i}n-1}\int_0^{t_1}s_1^{\frac2n}g_{k-i-1}(s_1,s) ds_1.
\end{align*}
Define
\begin{align*}
F_i(t_1,t_2,s)& = \frac n{2i} g_{k-i}(t_2,s) t_2^{\frac{2i}n} -\frac n{2i}g_{k-i}(t_1,s)t_1^{\frac{2i}n} + \frac{n^2}{2i(n-2i)}t_1^{\frac{2i-n}n}\int_0^{t_1}s_1^{\frac2n}g_{k-i-1}(s_1,s) ds_1\\
&\quad\quad -\frac{n^2}{2i(n-2i)}t_2^{\frac{2i-n}n}\int_0^{t_2}s_1^{\frac2n}g_{k-i-1}(s_1,s) ds_1
\end{align*}
we obtain \eqref{eq:recurrentclaimforg}. It follows easily from Lemma \ref{boundsforg} that
\begin{equation}\label{eq:proofclaim1}
\left(\int_0^{\infty} g_{k-i}(t,s)^{\frac n{n-2k}} ds\right)^{\frac{n-2k}n} \leq C t^{-\frac{2i}n},\quad \forall\, t> 0.
\end{equation}
It also follows from Lemma \ref{boundsforg} that if $s >t$ then
\[
\int_0^{t}s_1^{\frac2n}g_{k-i-1}(s_1,s) ds_1 \leq Ct^{\frac2n +1} s^{-1+\frac{2(k-i-1)}n},
\]
and if $t > s$ then
\begin{align*}
\int_0^{t}s_1^{\frac2n}g_{k-i-1}(s_1,s) ds_1& = \int_0^{s}s_1^{\frac2n}g_{k-i-1}(s_1,s) ds_1+\int_s^{t}s_1^{\frac2n}g_{k-i-1}(s_1,s) ds_1\\
&\leq C\int_0^ss_1^{\frac2n} s^{-1+\frac{2(k-i-1)}n} ds_1 + C \int_s^{t}s^{\frac2n} s_1^{-1+ \frac{2(k-i-1)}n} ds_1\\
&\leq C t^{\frac{2(k-i)}n}.
\end{align*}
Combining these two estimates, we get
\begin{equation}\label{eq:proofclaim2}
\left(\int_0^\infty \left(\int_0^{t} s_1^{\frac2n}g_{k-i-1}(s_1,s) ds_1\right)^{\frac n{n-2k}} ds\right)^{\frac{n-2k}n} \leq C t^{1-\frac{2i}n},\quad \forall\, t>0.
\end{equation}
Our desired estimate \eqref{eq:boundnormclaim} for $i=1,2,\ldots, k-2$ now is an immediate consequence of the estimates \eqref{eq:proofclaim1} and \eqref{eq:proofclaim2}. It remains to check \eqref{eq:forg1}. Indeed, using again integration by parts, we have
\begin{align*}
\int_{t_1}^{t_2} g(t,s) t^{\frac{2(k-1)}n-1} dt&= \frac{n}{2(k-1)} \int_{t_1}^{t_2} g(t,s) d t^{\frac{2(k-1)}n -1}\\
&= \frac{n}{2(k-1)} g(t_2,s)t_2^{\frac{2(k-1)}n} -\frac{n}{2(k-1)} g(t_1,s)t_1^{\frac{2(k-1)}n}\\
&\quad\quad \quad\quad+ \frac{n}{2(k-1)} \int_{t_1}^{t_2} s^{\frac2n} t^{\frac{2(k-1)}n-2} \chi_{\{s < t\}} dt
\end{align*}
Define
\[
F_{k-1}(t_1,t_2,s) = \frac{n}{2(k-1)} g(t_2,s)t_2^{\frac{2(k-1)}n} -\frac{n}{2(k-1)} g(t_1,s)t_1^{\frac{2(k-1)}n},
\]
we obtain \eqref{eq:forg1}. It also follows from Lemma \ref{boundsforg} that
\[
\left(\int_0^\infty g(t,s)^{\frac n{n-2k}} ds\right)^{\frac{n-2k}n} \leq C t^{-\frac{2(k-1)}n},\quad\forall\, t> 0,
\] 
which then immediately implies \eqref{eq:boundnormclaim}. Our claim is completely proved.

It is obvious that
\begin{align}\label{eq:last2}
\int_0^\infty &f^{**}(s) \int_{t_1}^{t_2} s^{\frac2n}t^{\frac{2(k-1)}n-2} \chi_{\{s < t\}} dt ds\notag\\
& =\int_0^{t_2} f^{**}(s) \int_{t_1}^{t_2} s^{\frac2n}t^{\frac{2(k-1)}n-2} \chi_{\{s < t\}} dt ds\notag\\
&=\int_0^{t_1} f^{**}(s) \int_{t_1}^{t_2} s^{\frac2n}t^{\frac{2(k-1)}n-2} \chi_{\{s < t\}} dt ds +\int_{t_1}^{t_2} f^{**}(s) \int_{t_1}^{t_2} s^{\frac2n}t^{\frac{2(k-1)}n-2} \chi_{\{s < t\}} dt ds\notag\\
&\leq \frac n{n-2(k-1)} t_1^{-1+\frac{2k}n}\int_0^{t_1}f^{**}(s) ds + \frac n{n-2(k-1)}\int_{t_1}^{t_2} f^{**}(s) s^{-1+ \frac{2k}n} ds.
\end{align}

Define
\begin{align*}
F(t_1,t_2,s)&= \sum_{i=1}^{k-1} c_{n,i} F_i(t_1,t_2,s).%+ c_{n,k} F_{k-1}(t_1,t_2,s).
\end{align*}
We then have from \eqref{eq:11111111}, \eqref{eq:recurrentclaimforg}, \eqref{eq:forg1} and \eqref{eq:last2} that
\begin{align}\label{eq:22222222}
(n\omega_n^{\frac1n})^{2k}(u^{*}(t_1)-u^*(t_2))& \leq c_{n,k}\int_{t_1}^{t_2}\frac{f^{**}(s)}{ s^{1-\frac{2k}n}} ds + \int_{t_1}^{t_2} F(t_1,t_2,s) f^{**}(s) ds \notag\\
&\quad\quad  c_{n,k} t_1^{-1+\frac{2k}n}\int_0^{t_1}f^{**}(s) ds.
\end{align}
Because of our claim above, we have
\begin{equation}\label{eq:33333333}
\left(\int_0^\infty F(t_1,t_2,s)^{\frac n{n-2k}} ds\right)^{\frac{n-2k}n} \leq C,
\end{equation}
with $C$ depends only on $n,k$.

H\"older inequality and Proposition \ref{normrearrangement} together implies
\begin{equation}\label{eq:last3}
\int_0^{t_1} f^{**}(s) ds \leq t_1^{1-\frac{2k}n} \left(\int_0^\infty (f^{**}(s))^{\frac n{2k}} ds\right)^{\frac{2k}n} \leq C \|f\|_{\frac{n}{2k}} t_1^{1-\frac{2k}n}
\end{equation}
with $C$ depends only on $n,k$. Thanks to \eqref{eq:22222222}, \eqref{eq:33333333}, \eqref{eq:last3}, H\"older inequality and Proposition \ref{normrearrangement}, we obtain \eqref{eq:keyinproof}, and hence finish our proof of Proposition \ref{keypropo}.

Let $u \in C_0^\infty(\R^n)$, denote $f = (-\Delta +I)^k u \in C_0^\infty(\R^n)$ then we have $u = L_{2k}\star f$, where for $\al \geq 0$
\[
L_\al (x) = \frac1{(4 \pi)^{\al/2}} \frac 1{\Gamma(\al/2)} \int_0^\infty e^{-\frac{\pi |x|^2}\de} e^{-\frac{\de}{4\pi}} \de^{-\frac{n - \al}2} \frac{d\de}{\de},
\]
is the Bessel potential. It is a classical fact that $u\in W^{2k,\frac n{2k}}(\R^n)$ with $k < \frac n2$ if and only if there exists $f \in L^{\frac n{2k}}(\R^n)$ such that $u =L_{2k} \star f$. From the definition of Bessel functions, we can readily check that $\int_{\R^n} L_\al(x) dx =1$ which then implies
\begin{equation}\label{eq:firsteigenvalue}
\inf\limits_{u\in W^{2k,\frac n{2k}}(\R^n), \, u\not\equiv 0} \frac{\|u\|_{2k,\frac n{2k}}^{\frac n{2k}}}{\|u\|_{\frac n{2k}}^{\frac n{2k}}} =1.
\end{equation}
The following result is an analogue result of Proposition \ref{keypropo} for the polyharmonic operator $(-\Delta + I)^k$. This will be point out to be crucial in the proof of Theorem \ref{Maintheoremfull} below.

\begin{proposition}\label{keyfull}
Let $u\in C_0^\infty(\R^n)$ and let $k$ be a positive integer less than $\frac n2$, denote $f = (-\Delta +I)^k u$, then there exists a constant $C(n,k)$ depending only on $n,k$ such that
\begin{equation}\label{eq:keyfull}
u^*(t_1) - u^*(t_2) \leq \frac{c(n,k)}{(n\omega_n^{\frac1n})^{2k}} \int_{t_1}^{t_2} \frac{f^{**}(s)}{s^{1 -\frac{2k}n}} ds + C(n,k)\|f\|_{\frac n{2k}},
\end{equation}
for any $0< t_1 < t_2< \infty$.
\end{proposition} 
To prove Proposition \ref{keyfull}, we need the following Lemma which is an analogue of \eqref{eq:MSrearrangement} for the operator $-\Delta +I$. More precisely, we prove the following estimate
\begin{lemma}\label{MStype}
Let $u\in C_0^\infty(\R^n)$ and denote $f = (-\Delta +I)u$, then 
\begin{equation*}
u^*(t_1) -u^*(t_2) \leq \frac1{(n\om_n^{1/n})^2} \int_{t_1}^{t_2} \frac{f^{**}(s)}{s^{1-\frac 2n}} ds,
\end{equation*}
for any $0 < t_1 < t_2 < \infty$.
\end{lemma}
With Lemma \ref{MStype} in hand, we can repeat the arguments in the proof of Proposition \ref{keypropo} to prove Proposition \ref{keyfull}. Hence we only have to prove Lemma \ref{MStype}. Our proof below follows the same lines in the proof of \eqref{eq:MSrearrangement} given in \cite{MS2014}, then we only sketch the proof and point out the differences.

\emph{Proof of Lemma \ref{MStype}:} For any $t, h > 0$, considering the test function $\phi$ defined in the proof of \cite[Lemma 3.3]{MS2014}, we have
\begin{align*}
\int_{\{t < |u| < t+h\}} |\na u|^2 dx & = \int_{\{t < |u| < t+h\}} (f-u)(|u|-t) \sign(u) dx +\int_{\{t+h < |u|\}} (f-u) h\sign(u) dx\\
&=\int_{\{t < |u| < t+h\}} f(|u|-t) \sign(u) dx +\int_{\{t+h < |u|\}} f h\sign(u) dx\\
&\quad \quad -\int_{\{t < |u| < t+h\}} |u|(|u|-t) dx -h \int_{\{t+h < |u|\}} |u| dx\\
&\leq \int_{\{t < |u| < t+h\}} |f|(|u|-t) dx +\int_{\{t+h < |u|\}} |f| h dx\\
&=\int_{\{t < |u|\}} |f|(|u|-t) dx -\int_{\{t+h < |u|\}} |f|(|u|-t-h) dx.
\end{align*}
Dividing both side of the latter inequality by $h$, letting $h \downarrow 0^+$, and then using Hardy-Littlewood inequality, we obtain
\[
-\frac{d}{dt}\int_{\{t < |u|\}} |\na u|^2 dx \leq \int_0^{\mu_u(t)} f^*(s) ds,\quad\forall\, t>0.
\]
Thanks to the latter estimate and \cite[Lemma 3.2]{MS2014}, by repeating the argument in the rest of the proof of \eqref{eq:MSrearrangement} in \cite{MS2014}, we get the desired result in Lemma \ref{MStype}.

Let us conclude this section by mentioning here that Proposition \ref{keypropo} incidentally gives a generalization of the Adachi--Tanaka inequality \eqref{eq:AdachiTanaka} to the space $W^{m,\frac nm}(\R^n)$ for any integer $m\in [1,n)$. This can be seen as the subcritical Adams inequality in $W^{m,\frac nm}(\R^n)$. We can state it as follows: for any $\al \in (0,\beta(n,m))$, there exists a constant $C(n,m,\al)$ depending only on $n,m$ and $\al$ such that
\begin{equation}\label{eq:AThigher}
\sup\limits_{u\in W^{m,\frac nm}(\R^n), \|\nabla^m u\|_{\frac nm} \leq 1} \int_{\R^n} \Phi_{n,m}(\al |u|^{\frac n{n-m}}) dx  \leq C(n,m,\al) \|u\|_{\frac nm}^{\frac nm}.
\end{equation}
This inequality recently was proved in \cite{FM2015}. Moreover, an upper bound for the constant $C(n,m,\al)$ can be found in that paper, i.e, $C(n,m,\al) \leq C/(\beta_{n,m} -\al)$, with $C$ depends only on $n$ and $m$. We can readily show that the inequality \eqref{eq:AThigher} with the previous upper bounds for $C(n,m,\al)$ is equivalent to the sharp Adams inequality \eqref{eq:FontanaMorpurgo}. We refer the reader to the paper \cite{LamLuZhang2015} for this interesting observation in the case $m=1,2$. Let us go to the proof of \eqref{eq:AThigher}. It is enough to prove it for compactly supported smooth functions. We first consider the case $m$ even, i.e, $m =2k$ for some $k\geq 1$. Suppose $u\in C_0^\infty(\R^n)$, let $f = \Delta^k u$. Choose $s_0\in (0,\infty)$ such that $u^*(s_0) =1$, then $s_0 \leq \|u\|_{\frac n{2k}}^{\frac n{2k}}$. For $s > s_0$ then $u^*(s) \leq 1$ which then implies 
\[
\Phi_{n,2k}(\al |u^*(s)|^{\frac n{n-2k}}) \leq C(n,2k,\al) |u^*(s)|^{\frac n{2k}}.
\]
Therefore we get
\[
\int_{s_0}^{\infty}\Phi_{n,2k}(\al |u^*(s)|^{\frac n{n-2k}}) ds \leq C(n,2k,\al)\|u\|_{\frac n{2k}}^{\frac n{2k}}.
\]
For $s < s_0$ we then have from Proposition \ref{keypropo} that
\[
u^*(s) \leq \lt(\frac1{\be(n,2k)} \ln \lt(\frac {s_0}s\rt)\rt)^{\frac{n-2k}n} + C(n,k)
\]
with $C(n,k)$ depends only on $n,k$. Using the elementary inequality \eqref{eq:elementary} for $p =\frac n{n-2k}$, we have
\[
|u^*(s)|^{\frac n{n-2k}} \leq \frac{1+\ep}{\beta(n,2k)} \ln \lt(\frac {s_0}s\rt) + C_\ep C(n,k)^{\frac n{n-2k}}.
\]
Choose $\ep >0$ small enough such that $\al(1+\ep) < \beta(n,2k)$, then we easily see that 
\[
\int_{0}^{s_0}\Phi_{n,2k}(\al |u^*(s)|^{\frac n{n-2k}}) ds  \leq C(n,2k,\al) s_0 \leq C(n,2k,\al)\|u\|_{\frac n{2k}}^{\frac n{2k}}.
\]
This finishes the proof when $m$ is even.

We next consider the case $m$ odd,i.e, $m=2k+1$ for some $k\geq 0$. It remains to prove for $k\geq 1$. Let $f= \Delta^k u$, then $\|f\|_{\frac n{2k}} \leq C$ with $C$ depends only on $n, k$by Sobolev inequality. Using \eqref{eq:welldone} below and the previous arguments we obtain \eqref{eq:AThigher} in this case.

\section{Proof of Theorem \ref{Maintheorem}}
With Proposition \ref{keypropo} in hand, we can follow the strategy of ${\rm\check{C}}$erny, Cianchi and Hencl \cite{Cerny2013} to prove Theorem \ref{Maintheorem}. We first reduce the proof of Theorem \ref{Maintheorem} to compactly supported smooth functions. This will be done by two following lemmas. The first one asserts that any function in $W^{m,\frac nm}_0(\Om)$ will satisfy an exponential integrability property.
\begin{lemma}\label{exponentialintegrability}
Let $u\in W^{m,\frac nm}_0(\Om)$ then for any $ p > 0$, it holds
\begin{equation}\label{eq:Orlicz}
\int_\Om e^{p |u|^{\frac n{n-m}}} dx < \infty.
\end{equation}
\end{lemma}
\begin{proof}
Choose a function $v\in C_0^\infty(\Om)$ such that $\|\nabla^m(u-v)\|_{\frac nm}^{\frac n{n-m}} p < \frac{\beta(n,m)}2$. Using the elementary inequality \eqref{eq:elementary} for $p= \frac n{n-m}$, we have
\[
|u|^{\frac n{n-m}} \leq 2 |u-v|^{\frac n{n-m}} + C_1 |v|^{\frac{n}{n-m}},
\]
hence by Adams inequality and the boundedness of function $v$, we obtain \eqref{eq:Orlicz}.
% From \eqref{eq:boundednearinfty} and \eqref{eq:Orlicz} we obtain Theorem \ref{Maintheorem} for general functions $u_j$.
\end{proof}

The second lemma enabes us reducing our proof to compactly supported smooth functions in $C_0^\infty(\Om)$. Its content is as follows.
\begin{lemma}\label{reducebounded}
Let $\{u_j\}_j$, $u$ and $p_1$ be as in Theorem \ref{Maintheorem}. Suppose that $\{v_j\}_j \in C_0^\infty(\Om)$ such that $\|\nabla^m (u_j -v_j)\|_{\frac nm} < \frac 1j$ for any $j$. Then for any $p_2 \in (p_1,P_{n,m}(u))$ there exists a constant $C$ which is independent of $j$ such that
\begin{equation*}
\sup_{j\geq 1}\int_\Om e^{p_1\beta(n,m) |u_j|^{\frac n{n-m}}} dx \leq C \sup_{j\geq 1} \int_\Om e^{p_2 \beta(n,m)|\tilde{v}_j|^{\frac n{n-m}}} dx,
\end{equation*}
where $\tilde{v}_j = \|\nabla^m u_j\|_{\frac nm} \frac{v_j}{\|\nabla^m v_j\|_{\frac nm}} \in C_0^\infty(\Om)$.
\end{lemma}
\begin{proof}
It is easy to see that $\|\nabla^m(u_j -\tilde{v}_j)\|_{\frac nm} \leq \frac 2j$ for any $j\geq 1$. It then implies that $v_j$ converges weakly to $u$ in $W^{m,\frac nm}_0(\Om)$. Using the elementary inequality \eqref{eq:elementary} for $p = \frac n{n-m}$ and $\ep = \frac{p_2-p_1}2$, we have
\[
|u_j|^{\frac n{n-m}} \leq (1+\ep) |\tilde{v}_j|^{\frac n{n-m}} + C_\ep |u_j -\tilde{v}_j|^{\frac n{n-m}}.
\]
Using H\"older inequality for $r = \frac{2p_2}{p_2+p_1} >1$, we obtain
\begin{align*}
\int_\Om e^{p_1\beta(n,m) |u_j|^{\frac n{n-m}}} dx &\leq \int_\Om e^{(1+\ep)p_1 \beta(n,m)|\tilde{v}_j|^{\frac n{n-m}} +p_1C_\ep \beta(n,m) |u_j-\tilde{v}_j|^{\frac n{n-m}}} dx\\
&\leq \lt(\int_\Om e^{p_2 \beta(n,m)|\tilde{v}_j|^{\frac n{n-m}}} dx\rt)^{\frac{p_2+p_1}{2p_2}} \lt(\int_\Om e^{\frac{r}{r-1}p_1C_\ep \beta(n,m) |u_j-\tilde{v}_j|^{\frac n{n-m}}} dx\rt)^{\frac{r-1} r}.
\end{align*}
Choose $J_0\in \N$ such that $J_0 \geq ((r p_1 C_\ep)/2(r-1))^{\frac {n-m}n}$, then for any $j \geq J_0$, by Adams inequality, we have
\[
\int_\Om e^{\frac{r}{r-1}p_1C_\ep \beta(n,m) |u_j-\tilde{v}_j|^{\frac n{n-m}}} dx \leq C(n,m) |\Om|
\]
with $C(n,m)$ depends only on $n,m$. Hence
\begin{align}\label{eq:boundednearinfty}
\sup_{j\geq J_0}\int_\Om e^{p_1\beta(n,m) |u_j|^{\frac n{n-m}}} dx &\leq (C(n,m) |\Om|)^{\frac{r-1}r} \lt(\sup_{j\geq J_0} \int_\Om e^{p_2 \beta(n,m)|\tilde{v}_j|^{\frac n{n-m}}} dx\rt)^{\frac{p_2+p_1}{2p_2}}\notag\\
&\leq C \sup_{j\geq J_0} \int_\Om e^{p_2 \beta(n,m)|\tilde{v}_j|^{\frac n{n-m}}} dx. 
\end{align}
%if Theorem \ref{Maintheorem} holds for compactly supported smooth functions in $\Om$.
Combining \eqref{eq:boundednearinfty} and Lemma \ref{exponentialintegrability} implies our desired estimate.
\end{proof}

%Indeed, for each $j\in \N$, we choose $v_j\in C_0^\infty(\Om)$ such that $\lt(\int_\Om |\nabla^m (u_j-v_j)|^{\frac nm} dx\rt)^{\frac mn} \leq \frac 1j$, then $v_j$ converges weakly to $u$ in $W_0^{m,\frac nm}(\Om)$.  
%In the other hand, for any function $u \in W_0^{m,\frac nm}(\Om)$, and for any $ p > 0$, it holds

%The rest of this section is devoted to prove Theorem \ref{Maintheorem} for compactly supported functions in $W^{m,\frac nm}_0(\Om)$. Since there are different arguments in the case $m$ even and $m$ odd, then we divide our proof into two cases following $m$ is even or odd. 
%\subsection{Proof of Theorem \ref{Maintheorem} for compactly supported smooths functions}
We continue our proof of Theorem \ref{Maintheorem} by contradiction argument. Because of Lemma \ref{reducebounded}. It is enough to prove Theorem \ref{Maintheorem} for compactly supported smooth functions. Suppose that there exists $\{u_j\}_j \subset C_0^\infty(\Om)$ such that $\|\nabla^m u_j\|_{\frac nm} \leq 1$ for any $j$, $u_j$ converges weakly to a non-zero function $u$ in $W^{m,\frac nm}_0(\Om)$, and $p_1 \in (1, P_{n,m}(u))$ such that
\begin{equation}\label{eq:contraassumpbounded}
\lim_{j\to \infty} \int_{\Om} e^{p_1 \beta(n,m) |u_j|^{\frac n{n-m}}} dx =\infty.
\end{equation}
Our aim is to look for a contradiction.

Using Rellich-Kondrachov theorem, we can additionally assume, up to a subsequence if necessary, that $u_j$ converges to $u$ almost everywhere in $\Om$ and in $L^p(\Om)$ for any $p < \infty$, and if $m$ is odd, we can make the same assumptions on functions $\Delta^{\frac{m-1}2} u_j$ and $\Delta^{\frac{m-1}2} u$ (but for $p < \frac{n}{m-1}$). We divide our proof into two cases.

{\bf Case $1$: Suppose that $m$ is even.} In this case, we can express $m=2k$ for some $k\geq 1$.
%We argue by contradiction. Suppose ther exist $p_1 \in (1,P_{n,2k}(u))$ and a sequence $\{u_j\}_j\in C_0^\infty(\Om)$ such that $\int_\Om |\Delta^k u_j|^{\frac n{2k}} dx \leq 1$ such that $u_j$ weakly converges to $u\not\equiv 0$ in $W_0^{2k,\frac n{2k}}(\Om)$ and 
%\begin{equation}\label{eq:contradictioneven}
%\lim_{j\to \infty} \int_\Om e^{p_1 \beta(n,2k) |u_j|^{\frac n{2k}}} dx =\infty.
%\end{equation} 
Denote $f_j = (-\Delta)^k u_j \in C_0^\infty(\Om)$ and $f=(-\Delta)^k u$, we then have $\|f_j\|_{\frac n{2k}} \leq 1$ for any $j\geq 1$, and $f_j$ converges weakly to $f$ in $L^{\frac n{2k}}(\Om)$. Lemma $2$ in \cite{doOMacedo2014} implies that, possibly passing to a subsequence, we then have $f_j^*$ converges a.e to a function $g$ in $(0, |\Om|)$, and
\begin{equation}\label{eq:normofg}
\int_0^{|\Om|} g(s)^{\frac n{2k}} ds \geq \int_0^{|\Om|} (f^*(s))^p ds = \int_\Om |f(x)|^p dx.
\end{equation}
Lemma \ref{keypropo} and integration by parts imply that
\begin{align*}
u_j^*(t_1) -u_j^*(t_2) &\leq \frac{c(n,k)}{(n\om_n^{1/n})^{2k}} \int_{t_1}^{t_2} \int_0^s f_j^*(t) dt s^{\frac{2k}n -2} ds + C(n,k)\\
&= \frac{n}{n-2k}\frac{c(n,k)}{(n\om_n^{1/n})^{2k}} \int_{t_1}^{t_2} \frac{f_j^*(s)}{s^{1-\frac{2k}n}} ds -\frac{n}{n-2k}\frac{c(n,k)}{(n\om_n^{1/n})^{2k}} t_2^{\frac{2k}n-1} \int_0^{t_2} f_j^*(s) ds \\
&\quad + \frac{n}{n-2k}\frac{c(n,k)}{(n\om_n^{1/n})^{2k}} t_1^{\frac{2k}n-1} \int_0^{t_1} f_j^*(s) ds + C(n,k)
\end{align*}
for any $0 < t_1 < t_2 < |\Om|$. H\"older inequality shows that
\[
t^{\frac{2k}n-1} \int_0^{t} f_j^*(s) ds \leq 1,\quad \forall\, t\in (0, |\Om|).
\]
We have proved the existence of a constant $C > 0$ depending only on $n,k$ such that 
\begin{equation}\label{eq:step1}
u_j^*(t_1) -u_j^*(t_2) \leq \frac{n}{n-2k}\frac{c(n,k)}{(n\om_n^{1/n})^{2k}} \int_{t_1}^{t_2} \frac{f_j^*(s)}{s^{1-\frac{2k}n}} ds + C,
\end{equation}
for any $0 < t_1 < t_2 < |\Om|$. Especially for $t_2 =|\Om|$, we obtain by using again H\"older inequality that
\[
u_j^*(t) \leq \lt(\frac1{\beta(n,2k)} \ln \lt(\frac{|\Om|}t\rt) \rt)^{1-\frac{2k}n} + C, \quad \forall\, t\in (0, |\Om|).
\]
%For any $a, b, \ep > 0$ and $p > 1$, we have the folowing elementary inequality 
%\begin{equation}\label{eq:elementary}
%(a+b) ^{\frac n{n-2k}} \leq (1+\ep) a^{\frac n{n-2k}} + C_\ep b^{\frac n{n-2k}}
%\end{equation}
Applying the elementary inequality \eqref{eq:elementary} for $ p = \frac n{n-2k}$, we obtain
\begin{equation}\label{eq:pointwiseboundeven}
(u_j^*(t))^{\frac n{n-2k}} \leq \frac{1+\ep}{\be(n,2k)} \ln \lt(\frac{|\Om|}t\rt) + C_\ep C^{\frac n{n-2k}},\quad \forall\, t\in (0, |\Om|).
\end{equation}
Note that \eqref{eq:pointwiseboundeven} yields the  Adams inequality \eqref{eq:Adams} with $\beta(n,2k)$ replaced by any smaller constant by choosing $\ep$ small enough.

We next claim that given any $p_2\in (p_1,P_{n,2k}(u))$, for every $j_0 \in \N$ and every $s_0\in (0, |\Om|)$, there is $j> j_0$ and $s\in (0,s_0)$ such that 
\begin{equation}\label{eq:observationeven}
u_j^*(s) \geq \lt(\frac{1}{p_2 \beta(n,2k)} \ln\lt(\frac{|\Om|}s\rt)\rt)^{\frac{n-2k}n}.
\end{equation}
Indeed, if this claim does not hold true, then there exist $j_0$ and $s_0$ such that 
\begin{equation}\label{eq:contra2even}
u_j^*(s) < \lt(\frac{1}{p_2 \beta(n,2k)} \ln\lt(\frac{|\Om|}s\rt)\rt)^{\frac{n-2k-1}n},
\end{equation}
for every $s\in (0, s_0)$ and $j\geq j_0$. It implies from \eqref{eq:contra2even} and \eqref{eq:pointwiseboundeven} that for any $j\geq j_0$ we have
\begin{align*}
\int_{\Om} e^{p_1 \be(n,2k) |u_j|^{\frac{n}{n-2k}}} dx &= \int_0^{|\Om|}e^{p_1 \be(n,2k) |u_j^*(s)|^{\frac{n}{n-2k}}} ds \\
&= \int_0^{s_0}e^{p_1 \be(n,2k) |u_j^*(s)|^{\frac{n}{n-2k}}} ds + \int_{s_0}^{|\Om|}e^{p_1 \be(n,2k) |u_j^*(s)|^{\frac{n}{n-2k}}} ds\\
&\leq \int_0^{s_0} \lt(\frac{|\Om|}s\rt)^{\frac{p_1}{p_2}} ds+e^{p_1 \beta(n,2k)C_\ep C_2^{\frac n{n-2k}}}\int_{s_0}^{|\Om|} \lt(\frac{|\Om|}s\rt)^{(1+\ep)p_1} ds.
\end{align*}
Therefore, we obtain
\[
\sup_{j} \int_{\Om} e^{p_1 \be(n,2k) |u_j|^{\frac{n}{n-2k}}} dx < \infty
\]
which contradicts with \eqref{eq:contraassumpbounded}. This proves our claim.

Thus, possibly passing to a subsequence, there exists a sequence $\{s_j\}$ such that 
\begin{equation}\label{eq:consequenceofclaimeven}
u_j^*(s_j) \geq \lt(\frac{1}{p_2 \beta(n,2k)} \ln\lt(\frac{|\Om|}{s_j}\rt)\rt)^{\frac{n-2k}n},\quad\text{and}\quad s_j \leq \frac 1j,\quad \forall\, j\in \N.
\end{equation}
Now, given $L > 0$ let us define the truncation operator $T^L$ and $T_L$ acting on any function $v$ by
\[
T^L(v) = \min \{|v|,L\} \sign(v),\quad \text{and}\quad  T_L(v) = v -T^L(v).
\]
It is easy to check that 
%\[
%\int_0^{|\Om|} |f_j^*(s)|^{\frac n{2k}} ds = \int_0^{|\Om|}|T^L(f_j^*)|^{\frac n{2k}} ds + \int_0^{|\Om|}|T_L(f_j^*)|^{\frac n{2k}} dx,
%\]
%and
\[
T^L(f_j^*)\longrightarrow T^L(g) \text{ a.e in } (0,|\Om|),\quad T_L(f_j^*)\longrightarrow T_L(g) \text{ a.e in } (0,|\Om|).
\]

For any $L > 0$, there exists a $j_0$ such that $u_j^*(s_j) > L$ for any $j \geq j_0$ because of \eqref{eq:consequenceofclaimeven}. This implies the existence of $r_j\in (s_j,|\Om|)$ such that $u_j^*(r_j) = L$. Moreover, it follows from \eqref{eq:step1} and \eqref{eq:consequenceofclaimeven} that there exists $C_1$ depending only on $n,k$ such that
\begin{align*}
\lt(\frac{1}{p_2 \beta(n,2k)} \ln\lt(\frac{|\Om|}{s_j}\rt)\rt)^{\frac{n-2k}n} - L & \leq u_j^*(s_j) -u_j^*(r_j)\\
&\leq \frac n{n-2k}\frac{c(n,k)}{(n\om_n^{1/n})^{2k}} \int_{s_j}^{r_j} \frac{f_j^*(s)}{s^{1-\frac{2k}n}} ds + C\\
&\leq \frac {n^2}{2k(n-2k)}\frac{c(n,k)}{(n\om_n^{1/n})^{2k}} f_j^*(s_j) (|\Om|^{\frac {2k}n} -s_j^{\frac{2k}n}) +C.
\end{align*}
This shows that $\lim_{j\to\infty} f_j^*(s_j) = \infty$. Hence, there exists $j_1 \geq j_0$ such that $f_j^*(s_j) > L$ for any $j\geq j_1$. Therefore, there exist $t_j \in (s_j,|\Om|)$ such that $f_j^*(t_j) = L$ and $f_j^*(s) < L$ for any $s > t_j$. Define $a_j =\min\{t_j, r_j\}$ for $j\geq j_1$, we have
\begin{align*}
&\lt(\frac{1}{p_2 \beta(n,2k)} \ln\lt(\frac{|\Om|}{s_j}\rt)\rt)^{\frac{n-2k}n} - L\\
&\quad\leq \frac n{n-2k}\frac{c(n,k)}{(n\om_n^{1/n})^{2k}} \int_{s_j}^{a_j} \frac{f_j^*(s)-L}{s^{1-\frac{2k}n}} ds + \frac n{n-2k}\frac{c(n,k)}{(n\om_n^{1/n})^{2k}} \int_{s_j}^{r_j} \frac{L}{s^{1-\frac{2k}n}} ds + C\\
&\quad \leq \frac n{n-2k}\frac{c(n,k)}{(n\om_n^{1/n})^{2k}} \int_{s_j}^{a_j} \frac{f_j^*(s)-L}{s^{1-\frac{2k}n}} ds + \frac {n^2}{2k(n-2k)}\frac{c(n,k)}{(n\om_n^{1/n})^{2k}}L (|\Om|^{\frac {2k}n} -s_j^{\frac{2k}n}) +C.
\end{align*}
Using H\"older inequality, we obtain
\begin{align*}
\frac n{n-2k}\frac{c(n,k)}{(n\om_n^{1/n})^{2k}}\int_{s_j}^{t_j} \frac{f_j^*(s)-L}{s^{1-\frac{2k}n}} ds&\leq \lt(\int_{s_j}^{t_j} \lt(f_j^*(s)-L\rt)^{\frac n{2k}} ds\rt)^{\frac{2k}{n}} \, \lt(\frac 1{\be(n,2k)} \ln \lt(\frac{t_j}{s_j}\rt)\rt)^{\frac{n-2k}n}\\
&\leq \lt(\int_{0}^{t_j} \lt(f_j^*(s)-L\rt)^{\frac n{2k}} ds\rt)^{\frac{2k}{n}} \, \lt(\frac 1{\be(n,2k)} \ln \lt(\frac{|\Om|}{s_j}\rt)\rt)^{\frac{n-2k}n}\\
&= \left(\int_0^{|\Om|}|T_L(f_j^*)|^{\frac n{2k}} ds\rt)^{\frac{2k}n} \lt(\frac 1{\be(n,2k)} \ln \lt(\frac{|\Om|}{s_j}\rt)\rt)^{\frac{n-2k}n}.
\end{align*}
For any $p_3 \in (p_2, P_{n,2k}(u))$, there exists $j_2 \geq j_1$ such that for any $j\geq j_2$
\[
p_3^{-\frac{n-2k}{2k}} \leq \int_0^{|\Om|}(T_L(f_j^*))^{\frac n{2k}} ds.
\]
Consequently, we have
\[
1- p_3^{-\frac{n-2k}{2k}} \geq \int_0^{|\Om|} \lt((f_j^*)^{\frac n{2k}} -(T_L(f_j^*))^{\frac n{2k}}\rt) ds, \quad\forall\, j\geq j_2.
\]
Since $f_j^* \to g$, and $T_L(f_j^*) \to T_L(g)$ a.e in $(0, |\Om|)$, respectively and $(f_j^*)^{\frac{n}{2k}} -(T_L(f_j^*))^{\frac{n}{2k}}$ is nonnegative functions, let $j$ tend to infinity and using Fatou lemma, we arrive
\[
1- p_3^{-\frac{n-2k}{2k}} \geq \int_0^{|\Om|} \lt(g^{\frac n{2k}} -(T_L(g))^{\frac n{2k}}\rt) ds,
\]
for any $L >0$ and . Letting $L$ tend to infinity, we obtain
\begin{equation*}%\label{eq:chanofp3}
1-\int_\Om |f|^{\frac n{2k}} dx < p_3^{-\frac{n-2k}{2k}} \leq 1- \int_0^{|\Om|} g^{\frac n{2k}} ds \leq 1-\int_\Om |f|^{\frac n{2k}} dx,
\end{equation*}
for any $p_3 \in (p_2, P_{n,2k}(u))$. This is impossible.

{\bf Case $2$: Suppose that $m$ is odd.} In this case, we can express $m =2k+1$ for some $k \geq 0$. Since the case $m=1$ (or $k=0$) was proved in \cite{Cerny2013}, hence we only prove for the case $m > 1$, i.e, $k\geq 1$. 

%We argue by contradiction. Suppose that there exist $p_1 \in (1,P_{n,2k+1}(u))$ and a sequence $\{u_j\}_j \in C_0^\infty(\Om)$ such that $\int_\Om |\nabla \Delta^k u_j|^{\frac n{2k+1}} dx \leq 1$, $u_j$ weakly converges to $u\not\equiv 0$ in $W_0^{2k+1, \frac n{2k+1}}(\Om)$ and
%\begin{equation}\label{eq:contradictionassumption}
%\lim_{j\to\infty} \int_{\Om} e^{p_1 \beta(n,2k+1) |u_j|^{\frac{n}{n-2k-1}}} dx = \infty.
%\end{equation}
%Using Rellich-Kondrachov theorem, possibly passing to a subsequence, we can assume that $f_j: = (-\Delta)^k u_j$ converges a.e to $f= (-\Delta)^k u$.

Denote $f_j = \Delta^k u_j$ and $f = \Delta^k u$, then we have $f_j \in W_0^{1,\frac n{2k+1}}(\Om)$ and $\|\nabla f_j\|_{\frac n{2k+1}}\leq 1$. The Sobolev inequality implies that $\|f_j\|_{\frac n{2k}} \leq C$ with $C$ is independent of $j$. It follows from Proposition \ref{keypropo} that
\begin{equation}\label{eq:applyofkeypropo}
u_j^*(t_1) -u_j^*(t_2) \leq \frac{c(n,k)}{(n\omega_n^{\frac1n})^{2k}} \int_{t_1}^{t_2} \frac{f_j^{**}(s)}{s^{1 -\frac{2k}n}} ds + C_1,
\end{equation}
with $C_1$ is independent of $j$, and for any $0 < t_1 < t_2 < |\Om|$.

The P\'olya-Szeg\"o principle implies that $\int_{\Om^\sharp} |\nabla f_j^\sharp|^{\frac n{2k+1}} dx \leq \int_\Om |\nabla f_j|^{\frac n{2k+1}} dx \leq 1$. This is equivalent to
\begin{equation}\label{eq:gradientbound}
\lt(\int_0^{|\Om|} \lt(n\om_n^{\frac1n} (-f_j^*)'(s)\rt)^{\frac n{2k+1}} s^{\frac{n-1}{2k+1}} ds \rt)^{\frac{2k+1}n}\leq 1.
\end{equation}
Since $f_j^*(|\Om|) =0$ and $f_j^*$ is locally absolutely continuous, then
\[
f_j^*(s) = -\int_s^{|\Om|} (f_j^*)'(r) dr,
\]
for $s\in (0,|\Om|)$. H\"older inequality and \eqref{eq:gradientbound} show that
\begin{align*}
f_j^*(s)& \leq \lt(\int_s^{|\Om|} \lt(n\om_n^{\frac1n} (-f_j^*)'(s)\rt)^{\frac n{2k+1}} s^{\frac{n-1}{2k+1}} ds\rt)^{\frac{2k+1}n} \lt(\int_s^{|\Om|} (n\om_n^{1/n})^{-\frac n{n-2k-1}} s^{-\frac{n-1}{n-2k-1}}\rt)^{\frac{n-2k-1}n}\\
&\leq C s^{-\frac {2k}n},
\end{align*}
with $C$ depends only on $n,k$. This shows that $\lim_{s\to 0} s f_j^*(s) =0$. Using integration by parts, we obtain
\begin{align*}
f_j^{**}(s) =\frac1s \int_0^s f_j^*(r) dr= f_j^*(s) + \frac 1s\int_0^s (-f_j^*)'(r) r dr,
\end{align*}
for any $s\in (0, |\Om|)$. Using again integration by parts we have
\begin{align*}
\int_{t_1}^{t_2} \frac{f_j^{**}(s)}{s^{\frac{2k}n -1}} ds &\leq \int_{t_1}^{t_2} s^{\frac{2k}n-1} f_j^*(s) ds + \int_{t_1}^{t_2} s^{\frac{2k}n-2} \int_0^s (-f_j^*)'(r) r dr ds\\
&=\frac{n}{2k} t_2^{\frac{2k}n} f_j^*(t_2) -\frac{n}{2k} t_1^{\frac{2k}n} f_j^*(t_1)  -\frac{n}{n-2k} t_2^{\frac{2k}n-1} \int_0^{t_2} (-f_j^*)'(r) r dr\\
&\quad + \frac{n}{n-2k} t_1^{\frac{2k}n-1} \int_0^{t_1} (-f_j^*)'(r) r dr + \frac {n^2}{2k(n-2k)} \int_{t_1}^{t_2} (-f_j^*)'(r) r^{\frac{2k}n} dr.
\end{align*}
It is easy, by using H\"older inequality and \eqref{eq:gradientbound}, to show that for any $t \in (0,|\Om|)$
\[
t^{\frac{2k}n} \int_{t}^{|\Om|} (-f_j^*)'(r) dr \leq C,\quad t^{\frac{2k}n-1} \int_{t}^{|\Om|} (-f_j^*)'(r)r dr \leq C,
\]
with $C$ is independent of $j$. Consequently, we have
\begin{equation}\label{eq:welldone}
u_j^*(t_1) -u_j^*(t_2) \leq \frac{c(n,k+1)}{(n\om_n^{1/n})^{2k}} \int_{t_1}^{t_2} (-f_j^*)'(r) r^{\frac{2k}n} dr + C_2,
\end{equation}
for any $0 < t_1 < t_2 < |\Om|$, with $C_2$ independent of $j$. Especially when $t_2 =|\Om|$, by using again H\"older inequality and \eqref{eq:gradientbound}, we obtain
\begin{equation*}
u_j^*(t) \leq \frac{c(n,k+1)}{(n\om_n^{1/n})^{2k+1}} \lt(\ln \lt(\frac{|\Om|}t\rt)\rt)^{\frac{n-2k-1}n} + C_2 = \lt(\frac1{\beta(n,2k+1)} \ln \lt(\frac{|\Om|}t\rt)\rt)^{\frac{n-2k-1}n} + C_2.
\end{equation*}
Using the elementary inequality \eqref{eq:elementary} for $p=\frac{n}{n-2k-1}$, we get
\begin{equation}\label{eq:pointwisebound}
(u_j^*(t))^{\frac{n}{n-2k-1}} \leq \frac{1+\ep}{\be(n,2k+1)} \lt(\frac{|\Om|}t\rt) + C_\ep C_2^{\frac{n}{n-2k-1}}.
\end{equation}
Note that \eqref{eq:pointwisebound} yields the  Adams inequality \eqref{eq:Adams} with $\beta(n,2k+1)$ replaced by any smaller constant by choosing $\ep$ small enough.

As in the case $m$ even, we next claim that given any $p_2\in (p_1,P_{n,2k+1}(u))$, for every $j_0 \in \N$ and every $s_0\in (0, |\Om|)$, there is $j> j_0$ and $s\in (0,s_0)$ such that 
\begin{equation}\label{eq:observation}
u_j^*(s) \geq \lt(\frac{1}{p_2 \beta(n,2k+1)} \ln\lt(\frac{|\Om|}s\rt)\rt)^{\frac{n-2k-1}n}.
\end{equation}
Indeed, if this claim does not hold, by repeating the argument in the case $m$ even and using \eqref{eq:pointwisebound}, we will get a contradiction of \eqref{eq:contraassumpbounded}.
% exist $j_0$ and $s_0$ such that 
%\begin{equation}\label{eq:contra2}
%u_j^*(s) < \lt(\frac{1}{p_2 \beta(n,2k+1)} \ln\lt(\frac{|\Om|}s\rt)\rt)^{\frac{n-2k-1}n},
%\end{equation}
%for every $s\in (0, s_0)$ and $j\geq j_0$, then by \eqref{eq:contra2} and \eqref{eq:pointwisebound}, we have for $j\geq j_0$
%\begin{align*}
%\int_{\Om} e^{p_1 \be(n,2k+1) |u_j|^{\frac{n}{n-2k-1}}} dx &= \int_0^{|\Om|}e^{p_1 \be(n,2k+1) |u_j^*(s)|^{\frac{n}{n-2k-1}}} ds \\
%&= \int_0^{s_0}e^{p_1 \be(n,2k+1) |u_j^*(s)|^{\frac{n}{n-2k-1}}} ds + \int_{s_0}^{|\Om|}e^{p_1 \be(n,2k+1) |u_j^*(s)|^{\frac{n}{n-2k-1}}} ds\\
%&\leq \int_0^{s_0} \lt(\frac{|\Om|}s\rt)^{\frac{p_1}{p_2}} ds+e^{p_1 \beta(n,2k+1)C_\ep C_2^{\frac n{n-2k-1}}}\int_{s_0}^{|\Om|} \lt(\frac{|\Om|}s\rt)^{(1+\ep)p_1} ds.
%\end{align*}
%The latter inequality shows that
%\[
%\sup_{j} \int_{\Om} e^{p_1 \be(n,2k+1) |u_j|^{\frac{n}{n-2k-1}}} dx < \infty
%\]
%which contradicts with \eqref{eq:contradictionassumption}. This proves our claim.

Thus, possibly passing to a subsequence, there exists a sequence $\{s_j\}$ such that 
\begin{equation}\label{eq:consequenceofclaim}
u_j^*(s_j) \geq \lt(\frac{1}{p_2 \beta(n,2k+1)} \ln\lt(\frac{|\Om|}{s_j}\rt)\rt)^{\frac{n-2k-1}n},\quad\text{and}\quad s_j \leq \frac 1j,\quad \forall\, j\in \N.
\end{equation}
Given $L > 0$, let us again use the truncation operators $T^L$ and $T_L$ defined above. It is easy to check that 
\[
\int_\Om |\nabla f_j|^{\frac n{2k+1}} dx = \int_{\Om}|\nabla T^L(f_j)|^{\frac n{2k+1}} dx + \int_{\Om}|\nabla T_L(f_j)|^{\frac n{2k+1}} dx,
\]
and
\[
T^L(f_j)\longrightarrow T^L(f) \text{ a.e in } \Om,\quad T_L(f_j)\longrightarrow T_L(f) \text{ a.e in } \Om.
\]

We next claim that for any $L > 0$, $T^L(f_j)$ converges weakly to $T^L(f)$ in $W^{1,\frac n{2k+1}}_0(\Om)$. Indeed, $\{T^L(f_j)\}_j$ is bounded in $W^{1,\frac n{2k+1}}_0(\Om)$, hence any its subsequence possesses a subsequence which weakly converges to some function in $W^{1,\frac n{2k+1}}_0(\Om)$. Since this subsequence converges almost everywhere to $T^L(f)$ (and converges in $L^p(\Om)$ for any $p < \frac n{2k}$), so the weak limit function must be $T^L(f)$. This proves our claim. It is obvious that $T_L(f_j)$ also converges weakly to $T_L(g)$ in $W^{1,\frac n{2k+1}}_0(\Om)$.

For any $L >0$, there exists $j_0 >0$ such that $u_j^*(s_j) > L$ for every $j\geq j_0$ because of \eqref{eq:consequenceofclaim}. Consequently, there exists $r_j \in (s_j, |\Om|)$ such that $u_j^*(r_j) =L$. It follows from \eqref{eq:welldone} and \eqref{eq:consequenceofclaim} that
\begin{align*}
\lt(\frac{1}{p_2 \beta(n,2k+1)} \ln\lt(\frac{|\Om|}{s_j}\rt)\rt)^{\frac{n-2k-1}n} - L & \leq u_j^*(s_j) -u_j^*(r_j)\\
&\leq \frac{c(n,k+1)}{(n\om_n^{1/n})^{2k}} \int_{s_j}^{r_j} (-f_j^*)'(r) r^{\frac{2k}n} dr + C_2\\
&\leq \frac{c(n,k+1)}{(n\om_n^{1/n})^{2k}} |\Om|^{\frac {2k}n} f_j^*(s_j) + C_2,
\end{align*}
which implies $\lim_{j\to \infty}f_j^*(s_j) = \infty$. Hence there exists $j_1 \geq j_0$ such that $f_j^*(s_j) > L$ for any $j\geq j_1$, so there exists $t_j\in (s_j, |\Om|)$ such that $f_j^*(t_j) = L$. Denote $a_j =\min\{t_j,r_j\}$ we have
\begin{align*}
u_j^*(s_j) -u_j^*(r_j)&\leq \frac{c(n,k+1)}{(n\om_n^{1/n})^{2k}} \int_{s_j}^{a_j} (-f_j^*)'(r) r^{\frac{2k}n} dr+ \frac{c(n,k+1)}{(n\om_n^{1/n})^{2k}} \int_{a_j}^{r_j} (-f_j^*)'(r) r^{\frac{2k}n} dr + C_2\\
&\leq \frac{c(n,k+1)}{(n\om_n^{1/n})^{2k}} \int_{s_j}^{a_j} (-f_j^*)'(r) r^{\frac{2k}n} dr+ \frac{c(n,k+1)}{(n\om_n^{1/n})^{2k}} r_j^{\frac {2k}n} (f_j^*(a_j)-f_j^*(r_j)) + C_2\\
&\leq \frac{c(n,k+1)}{(n\om_n^{1/n})^{2k}} \int_{s_j}^{a_j} (-f_j^*)'(r) r^{\frac{2k}n} dr+ \frac{c(n,k+1)}{(n\om_n^{1/n})^{2k}} |\Om|^{\frac {2k}n} L + C_2,
\end{align*}
here we use the fact that if $t_j < r_j$ then $f_j^*(a_j)-f_j^*(r_j) \leq L$ while if $t_j\geq r_j$ then $f_j^*(a_j)-f_j^*(r_j) =0$. A simple application of H\"older inequality yields
\begin{align*}
\frac{c(n,k+1)}{(n\om_n^{1/n})^{2k}} &\int_{s_j}^{a_j} (-f_j^*)'(r) r^{\frac{2k}n} dr\\
&\leq \lt(\int_{s_j}^{a_j} \lt((-f_j^*)'(r) n\om_n^{\frac1n}\rt)^{\frac n{2k+1}} r^{\frac{n-1}{2k+1}} dr\rt)^{\frac{2k+1}{n}} \, \lt(\frac 1{\be(n,2k+1)} \ln \lt(\frac{a_j}{s_j}\rt)\rt)^{\frac{n-2k-1}n}\\
&\leq \lt(\int_{0}^{t_j} \lt((-f_j^*)'(r) n\om_n^{\frac1n}\rt)^{\frac n{2k+1}} r^{\frac{n-1}{2k+1}} dr\rt)^{\frac{2k+1}{n}} \, \lt(\frac 1{\be(n,2k+1)} \ln \lt(\frac{|\Om|}{s_j}\rt)\rt)^{\frac{n-2k-1}n}\\
&= \lt(\int_\Om |\nabla T_L(f_j^\sharp)|^{\frac n{2k+1}}dx\rt)^{\frac{2k+1}n} \lt(\frac 1{\be(n,2k+1)} \ln \lt(\frac{|\Om|}{s_j}\rt)\rt)^{\frac{n-2k-1}n}.
\end{align*}
Hence , for any $p_3\in (p_2, P_{n,2k+1}(u))$, there exists $j_2\geq j_1$ such that
\begin{equation}\label{eq:almostdone}
p_3^{-\frac{n-2k-1}n} \leq \lt(\int_\Om |\nabla T_L(f_j^\sharp)|^{\frac n{2k+1}}dx\rt)^{\frac{2k+1}n},\quad\forall\, j\geq j_2.
\end{equation}
Note that $T_L(f_j^\sharp) = (T_L(f_j))^\sharp$, then by P\'olya-Szeg\"o principle, we have
\[
\lt(\int_\Om |\nabla T_L(f_j)|^{\frac n{2k+1}}dx\rt)^{\frac{2k+1}n} \geq \lt(\int_\Om |\nabla (T_L(f_j))^\sharp|^{\frac n{2k+1}}dx\rt)^{\frac{2k+1}n} =\lt(\int_\Om |\nabla T_L(f_j^\sharp)|^{\frac n{2k+1}}dx\rt)^{\frac{2k+1}n}
\]
This inequality and \eqref{eq:almostdone} imply that
\[
p_3^{-\frac{n-2k-1}n} \leq \lt(\int_\Om |\nabla T_L(f_j)|^{\frac n{2k+1}}dx\rt)^{\frac{2k+1}n}, \quad\forall\, j\geq j_2
\]
or equivalently,
\[
p_3^{-\frac{n-2k-1}{2k+1}} \leq \int_\Om |\nabla T_L(f_j)|^{\frac n{2k+1}}dx \leq 1 -\int_\Om |\nabla T^L(f_j)|^{\frac n{2k+1}}dx,\quad\forall\, j\geq j_2.
\]
Therefore, we obtain
\begin{equation}\label{eq:almostdone1}
1- p_3^{-\frac{n-2k-1}{2k+1}} \geq \int_\Om |\nabla T^L(f_j)|^{\frac n{2k+1}}dx,\quad\forall\, j\geq j_2
\end{equation}
Letting $j\to \infty$ in \eqref{eq:almostdone1} and using the weak lower semicontinuity of $L^{\frac n{2k+1}}-$norm of gradient, we obtain
\begin{equation*}\label{eq:almostdone2}
1- p_3^{-\frac{n-2k-1}{2k+1}} \geq \int_\Om |\nabla T^L(f)|^{\frac n{2k+1}}dx,
\end{equation*}
for any $p_3 \in (p_2, P_{n,2k+1}(u)$ and any $L >0$. Letting $L$ to infinity, we get
\[
1- p_3^{-\frac{n-2k-1}{2k+1}} \geq \int_\Om |\nabla f|^{\frac n{2k+1}}dx,
\]
for any $p_3 \in (p_2, P_{n,2k+1}(u))$. This is impossible.

%Combining \eqref{eq:almostdone2} and \eqref{eq:Llarge}, we have
%\begin{align*}
%p_3 \geq \lt(1-\int_\Om |\nabla T^L(f)|^{\frac n{2k+1}}dx\rt)^{-\frac{2k+1}{n-2k-1}} > \frac{p_3}{P_{n,2k+1}(u)}\lt( \frac 1{1-\int_\Om |\nabla f|^{\frac n{2k+1}}dx}\rt)^{\frac{2k+1}{n-2k-1}} =p_3
%\end{align*}
%which is a contradiction.

%We have two following cases:

%{\bf Case $1$: $0 < \int_\Om |\nabla f|^{\frac n{2k+1}}dx < 1.$} Fix any $p_3 \in (p_2, P_{n,2k+1}(u)$, since
%\[
%\lim_{L\to\infty} \int_\Om |\nabla T^L(f)|^{\frac n{2k+1}}dx =\lim_{L\to\infty} \int_{\{|f|\leq L\}} |\nabla f|^{\frac n{2k+1}}dx =\int_\Om |\nabla f|^{\frac n{2k+1}}dx
%\]
%hence we can choose a $L > 0$ so large such that
%\begin{equation}\label{eq:Llarge}
%\frac{1 -\int_\Om |\nabla f|^{\frac n{2k+1}}dx}{1 -\int_\Om |\nabla T^L(f)|^{\frac n{2k+1}}dx} > \lt(\frac{p_3}{P_{n,2k+1}(u)}\rt)^{\frac{n-2k-1}{2k+1}}.
%\end{equation}

%{\bf Case $2$: $\int_\Om |\nabla f|^{\frac n{2k+1}}dx =1$.} The proof proceeds along the same %lines in the proof of {\bf Case $1$}, and we limit ourselves to sketching a few differences. We %fix any $p_3 > p_2$ and choose $L > 0$ large enough such that 
%\[
%\int_\Om |\nabla T^L(f)|^{\frac n{2k+1}}dx > 1 -\frac12 \lt(\frac 1{p_3}\rt)^{\frac{n-2k-1}{2k+1}}.
%\]
%For such a $L$, repeating the argument in {\bf Case $1$} we see that \eqref{eq:almostdone2} also holds. Therefore we get a contradiction that
%\[
%1 -\lt(\frac 1{p_3}\rt)^{\frac{n-2k-1}{2k+1}}\geq \int_\Om |\nabla T^L(f)|^{\frac n{2k+1}}dx > 1 -%\frac12 \lt(\frac 1{p_3}\rt)^{\frac{n-2k-1}{2k+1}}.
%\]
%Our proof hence is finished.
It remains to check the sharpness of the exponent $P_{n,m}(u)$. We will show that for any $\al \in (0,1)$, there exists a sequence $\{u_j\}_j\subset W^{m,\frac nm}_0(\Om)$ and $u\in W^{m,\frac nm}_0(\Om)$ such that $\|\na^m u_j\|_{\frac nm} =1$, $\|\na^m u\|_{\frac nm} =\al$, $u_j \rightharpoonup u$ in $W^{m,\frac nm}_0(\Om)$, $u_j\to u$ a.e in $\Om$ such that
\[
\lim_{j\to\infty}\int_\Om e^{\beta(n,m) (1-\al^{\frac nm})^{-\frac m{n-m}} |u_j|^{\frac n{n-m}}} dx =\infty.
\]
By scaling, we can assume that $\o{B}_2 \subset \Om$. For $j\geq 2$, define
\[
v_j(x) = 
\begin{cases}
\lt(\frac {\ln j}{\beta(n,m)}\rt)^{1-\frac mn} + \frac{n\beta(n,m)^{\frac mn -1}}{2 (\ln j)^{\frac mn}}\sum_{l=1}^{m-1} \frac{(1 -j^{\frac 2n} |x|^2)^l}{l}&\mbox{if $0\leq |x|\leq j^{-\frac 1n}$}\\
-n \beta(n,m)^{\frac mn -1} (\ln j)^{-\frac mn} \ln |x| &\mbox{if $j^{-\frac 1n} \leq |x| < 1$}\\
\xi_j(x) &\mbox{if $1\leq |x|\leq 2$},
\end{cases}
\]
where $\xi_j\in C_0^\infty(B_2)$ are radial functions which are chosen such that $\xi_j =0$ on $\partial B_1$ and $\partial B_2$, and for $l =1,2,\ldots,m-1$
\[
\frac{\partial^l \xi_j}{\partial r^l} \bigl{|}_{\partial B_1} = (-1)^l (l-1)!\beta(n,m)^{\frac mn -1} (\ln j)^{-\frac mn}, \quad \quad \frac{\partial^l \xi_j}{\partial r^l} \bigl{|}_{\partial B_2} =0,
\]
and $\xi_j$, $|\nabla^l \xi_j|$ and $|\nabla^m \xi_j|$ are all $O((\ln j)^{-\frac mn})$. The choice of these functions is inspired from \cite{Zhao2013}.

An easy computation shows that 
\[
1 \leq \|\nabla^m v_j\|_{\frac nm}^{\frac nm} \leq 1 + O((\ln j)^{-1}).
\]
Setting $\tilde{v}_j = v_j /\|\nabla^m v_j\|_{\frac nm}$, then we have $\tilde{v}_j \rightharpoonup 0$ in $W^{m,\frac nm}_0(\Om)$ and converges a.e to $0$ in $\Om$. Taking a function $v\in C_0^\infty(\Om)$ such $v$ is constant in $B_2$ and $\|\nabla^m v\|_{\frac nm} = \al$. Define $u_j = v + (1-\al^{\frac nm})^{\frac mn} \tilde{v_j}$ then $u_j\in W_0^{m,\frac nm}(\Om)$, $\|\nabla^m u_j\|_{\frac nm} = 1$ for all $j\geq 2$ since the supports of $\nabla^m v$ and $\nabla^m \tilde{v}_j$ are disjoint, and $u_j \rightharpoonup v$ in $W^{m,\frac nm}_0(\Om)$. Replacing $v$ by $-v$ if necessary, we can assume that $v\geq A$ on $B_2$ for some $A >0$. Therefore we have
\begin{align*}
&\int_{\Om} e^{\beta(n,m) (1-\al^{\frac nm})^{-\frac m{n-m}} |u_j|^{\frac n{n-m}}} dx\\
&\geq \int\limits_{|x| \leq j^{-\frac 1n}} \exp\lt(\frac{\beta(n,m)}{ (1-\al^{\frac nm})^{\frac m{n-m}}} \lt(A +\frac{(1-\al^{\frac nm})^{\frac mn}}{(1 + O((\ln j)^{-1}))^{\frac mn}}\lt(\frac {\ln j}{\beta(n,m)}\rt)^{1-\frac mn}\rt)^{\frac n{n-m}}\rt) dx\\
&\geq \om_n \exp\lt(\lt(C + \frac{(\ln j)^{\frac{n-m}n}}{(1 + O((\ln j)^{-1}))^{\frac mn}}\rt)^{\frac{n}{n-m}} -\ln j\rt), 
\end{align*} 
for some $C>0$ which is independent of $j$. It is easy to see that there exists a constant $0 < C_1 < C$ and $j_0$ such that for any $j\geq j_0$, we have
\[
\exp\lt(\lt(C + \frac{(\ln j)^{\frac{n-m}n}}{(1 + O((\ln j)^{-1}))^{\frac mn}}\rt)^{\frac{n}{n-m}} -\ln j\rt) \geq \exp\lt(\lt(C_1 + (\ln j)^{\frac{n-m}n}\rt)^{\frac{n}{n-m}} -\ln j\rt),
\]
hence
\[
\liminf_{j\to\infty} \int_{\Om} e^{\beta(n,m) (1-\al^{\frac nm})^{-\frac m{n-m}} |u_j|^{\frac n{n-m}}} dx \geq \lim_{j\to\infty} \exp\lt(\lt(C_1 + (\ln j)^{\frac{n-m}n}\rt)^{\frac{n}{n-m}} -\ln j\rt) = \infty.
\]

\section{Proof of Theorem \ref{Maintheoremfull} and Theorem \ref{Maintheorem2}}
\subsection{Proof of Theorem \ref{Maintheoremfull}}
We use again the contradiction argument to prove Theorem \ref{Maintheoremfull}. We first show that we can reduce the proof to the compactly supported smooth functions in $\R^n$. This step is done by two following lemmas.
\begin{lemma}\label{expointegral}
Let $u\in W^{m,\frac nm}(\R^n)$ then for any $p >0$, it holds
\[
\int_{\R^n} \Phi_{n,m}(p |u|^{\frac n{n-m}}) dx < \infty.
\]
\end{lemma}
\begin{proof}
For $\ep >0$, we choose $v\in C_0^\infty(\R^n)$ such that $\|u-v\|_{m,\frac nm} < \ep$. For any $\de >0$, applying the elementary inequality \eqref{eq:elementary} for $p = \frac n{n-m}$ we have
\[
|u|^{\frac n{n-m}} \leq (1+\de) |u-v|^{\frac n{n-m}} + C_{\de} |v|^{\frac{n}{n-m}}.
\]
Let us divide $\R^n$ into two parts as follows
\begin{gather*}
\Om_1 =\{x\, :\, |u(x)-v(x)| \leq 1\},\, \Om_2 =\{x\, :\, |u(x)-v(x)| > 1\},
\end{gather*}

On $\Om_1$, we have $|u| \leq 1+\max\{|v(x)|\, :\, x\in \R^n\} =:C_v$ then $\Phi_{n,m}(p |u|^{\frac{n}{n-m}}) \leq C(n,m,p,v)|u|^{\frac nm}$ for some $C(n,m,p,v)>0$ depending on $n, m,p$ and $v$, hence
\begin{align*}
\int_{\Om_1} \Phi_{n,m}(p |u|^{\frac{n}{n-m}}) dx &\leq \int_{\{|u|\leq C_v\}}\Phi_{n,m}(p |u|^{\frac{n}{n-m}}) dx\\
&\leq C(n,m,p,v)\int_{\{|u|\leq C_v\}} |u|^{\frac nm} dx\\
&\leq C(n,m,p,v) C \|u\|_{m,\frac nm}^{\frac nm}\\
&< \infty,
\end{align*}
here we use the inequality \eqref{eq:firsteigenvalue}. 

%It is easy to prove that for any $a, b> 0$ there exist a constant $C(n,m,a,b) > 0$ such that 
%\begin{equation}\label{eq:stepineq}
%\Phi_{n,m}(s+t) \leq C(n,m,a,b) \Phi_{n,m}(t),\quad \forall\, s\in [0,a],\, \forall\, t\geq b.
%\end{equation}
%Since $v$ is bounded, hence there exists a constant $C(n,m,p,\de,v)$ depending on $n,m,p,\de$ and $v$ such that for any $x\in \Om_2$
%\[
%\Phi_{n,m}(p(1+\de) |u-v|^{\frac{n}{n-m}} +pC_\de|v|^{\frac n{n-m}})\leq C(n,m,p,\de,v) \Phi_{n,m}(p(1+\de) |u-v|^{\frac{n}{n-m}}),
%\]
%which then implies
For the integral on $\Om_2$, using the elementary inequality \eqref{eq:elementary} for $\ep =1$ and $p =\frac n{n-m}$, we have
\begin{align*}
\int_{\Om_2}\Phi_{n,m}(p |u|^{\frac{n}{n-m}}) dx &\leq \int_{\Om_2} \exp(p|u|^{\frac{n}{n-m}}) dx\\
&\leq  \int_{\Om_2} \exp(2p |u-v|^{\frac{n}{n-m}} +pC_1|v|^{\frac n{n-m}}) dx\\
&\leq C(n,m,p,v) \int_{\R^n} \Phi_{n,m}(2p |u-v|^{\frac{n}{n-m}}) dx,
\end{align*}
here we use the fact $|u-v| \geq 1$ on $\Om_2$ and $v$ is bounded. Choosing $\ep$ small enough such that $2p\ep^{\frac{n}{n-m}} \leq \beta(n,m)$ and using the Adams inequality \eqref{eq:FontanaMorpurgo}, we get
\begin{align*}
\int_{\Om_2} \Phi_{n,m}(2p |u-v|^{\frac{n}{n-m}}) dx\leq \int_{\R^n} \Phi_{n,m}\lt(\beta(n,m) \lt(\frac{|u-v|}{\ep}\rt)^{\frac n{n-m}}\rt)dx
< \infty,
\end{align*}
which then implies
\[
\int_{\Om_2}\Phi_{n,m}(p |u|^{\frac{n}{n-m}}) dx < \infty.
\]
The proof of Lemma \ref{expointegral} hence is finished.
%We have shown that
%\[
%\int_{\R^n}\Phi_{n,m}(p |u|^{\frac{n}{n-m}}) dx =\int_{\Om_1}\Phi_{n,m}(p |u|^{\frac{n}{n-m}}) dx +\int_{\Om_2}\Phi_{n,m}(p |u|^{\frac{n}{n-m}}) dx< \infty.
%\]
\end{proof}
\begin{lemma}\label{reduceunbounded}
Let $\{u_j\}_j$, $u$ and $p_1$ be as in Theorem \ref{Maintheoremfull}. Suppose that $\{v_j\}_j \in C_0^\infty(\R^n)$ such that $\|u_j -v_j\|_{m,\frac nm} < \frac 1j$ for any $j$. Then for any $p_2 \in (p_1,P_{n,m}(u))$ there exists a constant $C$ which is independent of $j$ such that
\begin{equation*}
\sup_{j\geq 1}\int_{\R^n}\Phi_{n,m}(p_1\beta(n,m) |u_j|^{\frac n{n-m}}) dx \leq C \sup_{j\geq 1} \int_{\R^n}\Phi_{n,m}(p_2 \beta(n,m)|\tilde{v}_j|^{\frac n{n-m}}) dx,
\end{equation*}
where $\tilde{v}_j = \| u_j\|_{m,\frac nm} \frac{v_j}{\|v_j\|_{m,\frac nm}} \in C_0^\infty(\R^n)$.
\end{lemma}
\begin{proof}
¨%It is easy to see that for any $A > 0$, there is a constant $C(n,m,A)$ depending only on $n,m,A$ such that $\phi_{n,m}(t) \leq C(n,m,A) t^{\frac {n-m}m}$ for any $t\leq A$. 
It is easy to see that
\begin{equation}\label{eq:mienkobichan}
\int_{\{|u_j| \leq 2\}} \Phi_{n,m}(p_1 \beta(n,m) |u_j|^{\frac n{n-m}}) dx \leq C \int_{\R^n} |u_j|^{\frac nm} dx \leq C,
\end{equation}
where $C$ is independent of $j$, here we use again the inequality \eqref{eq:firsteigenvalue}. We next divide the set $\{|u_j| > 2\}$ into two parts as follows
\[
\Om_{j,1} =\{|u|_j >2\} \cap \{|u_j -\tilde{v}_j| \leq 1\},\quad \Om_{j,2} = \{|u|_j >2\} \cap \{|u_j -\tilde{v}_j| > 1\}.
\]
On $\Om_{j,1}$ we have $|\tilde{v}_j| \geq  |u_j| -|u_j-\tilde{v}_j| > 1$. Then using elementary inequality \eqref{eq:elementary} for $\de = \frac{p_2}{p_1} -1$, we obtain
\begin{align}\label{eq:Omegaj1}
\int_{\Om_{j,1}} \Phi_{n,m}(p_1 \beta(n,m) |u_j|^{\frac n{n-m}}) dx &\leq\int_{\Om_{j,1}} \exp(p_1 \beta(n,m) |u_j|^{\frac n{n-m}}) dx\notag\\
&\leq \int_{\Om_{j,1}} \exp(p_1 \beta(n,m) (1+\de) |\tilde{v}_j|^{\frac n{n-m}}+ p_1\beta(n,m) C_\de) dx\notag\\
&\leq C \int_{\R^n} \Phi_{n,m}(p_2 \beta(n,m) |\tilde{v}_j|^{\frac n{n-m}}) dx\notag\\
&\leq C \sup_{j \in \N} \int_{\R^n} \Phi_{n,m}(p_2 \beta(n,m) |\tilde{v}_j|^{\frac n{n-m}}) dx,
\end{align}
here we use the fact that $|\tilde{v}_j| \geq 1 $ on $\Om_{j,1}$. 

We next estimate the integral on $\Om_{j,2}$. To do this, we divide it into two parts as follows
\[
\Om_{j,2}^1 =\Om_{j,2} \cap \{|\tilde{v}_j| < 1\},\quad \Om_{j,2}^2 = \Om_{j,2}^2 = \Om_{j,2} \cap \{|\tilde{v}_j| \geq 1\}.
\]
On $\Om_{j,2}^1$, using the elementary inequality \eqref{eq:elementary} for $\ep =1$, we have
\begin{align*}
\int_{\Om_{j,2}^1} \Phi_{n,m}(p_1 \beta(n,m) |u_j|^{\frac n{n-m}}) dx &\leq\int_{\Om_{j,2}^1} \exp(p_1 \beta(n,m) |u_j|^{\frac n{n-m}}) dx\notag\\
&\leq \int_{\Om_{j,2}^2} \exp(2p_1 \beta(n,m) |u_j-\tilde{v}_j|^{\frac n{n-m}}+ p_1\beta(n,m) C_1) dx\notag\\
&\leq C \int_{\R^n} \Phi_{n,m}(2p_1 \beta(n,m)|u_j-\tilde{v}_j|^{\frac n{n-m}}) dx.
\end{align*}
It is obvious that $\|u_j -\tilde{v}_j\|_{m,\frac nm} \leq \frac 2j$ for any $j$. Choose $J_0$ such that $J_0 \geq (4p)^{\frac{n-m}n}$, we have for any $j \geq J_0$ that
\begin{align*}
\int_{\R^n} \Phi_{n,m}(2p_1 \beta(n,m)|u_j-\tilde{v}_j|^{\frac n{n-m}}) dx &\leq \int_{\R^n} \Phi_{n,m}\lt( \beta(n,m)\lt|\frac{u_j-\tilde{v}_j}{\|u_j-\tilde{v}_j\|_{m,\frac nm}}\rt|^{\frac n{n-m}}\rt) dx\\
&\leq C.
\end{align*}
These estimates and Lemma \ref{expointegral} imply that
\begin{equation}\label{eq:Omegaj21}
\sup_{j\in \N}\int_{\Om_{j,2}^1} \Phi_{n,m}(p_1 \beta(n,m) |u_j|^{\frac n{n-m}}) dx \leq C
\end{equation}
On $\Om_{j,2}^2$, using the elementary inequality \eqref{eq:elementary} for $\ep = \frac{p_2-p_1}{2p_1}$, we have 
\[
|u_j|^{\frac n{n-m}} \leq (1+\ep) |\tilde{v}_j|^{\frac n{n-m}} + C_\ep |u_j-\tilde{v}_j|^{\frac n{n-m}}.
\]
Denote $r = \frac{2p_2}{p_1+p_2}$ and $r' = \frac{r}{r-1}$, using H\"older inequality, we have
\begin{align*}
&\int_{\Om_{j,2}^2} \Phi_{n,m}(p_1 \beta(n,m) |u_j|^{\frac n{n-m}}) dx \leq\int_{\Om_{j,2}^2} \exp(p_1 \beta(n,m) |u_j|^{\frac n{n-m}}) dx\notag\\
&\leq \int_{\Om_{j,2}^2} \exp((1+\ep)p_1 \beta(n,m) |\tilde{v}_j|^{\frac n{n-m}}+ p_1\beta(n,m) C_\ep |u_j-\tilde{v}_j|^{\frac{n}{n-m}}) dx\notag\\
&\leq \lt(\int_{\Om_{j,2}^2} \exp(p_2 \beta(n,m)|\tilde{v}_j|^{\frac n{n-m}}) dx\rt)^{\frac{p_2+p_1}{2p_2}} \lt(\int_{\Om_{j,2}^2} \exp(r'p_1C_\ep \beta(n,m) |u_j-\tilde{v}_j|^{\frac n{n-m}})dx\rt)^{\frac 1{r'}}\\
&\leq C \lt(\int_{\R^n} \Phi_{n,m}(p_2 \beta(n,m)|\tilde{v}_j|^{\frac n{n-m}}) dx\rt)^{\frac{p_1+p}{2p_1}} \lt(\int_{\R^n} \Phi_{n,m}(r'p_1C_\ep \beta(n,m) |u_j-\tilde{v}_j|^{\frac n{n-m}})dx\rt)^{\frac 1{r'}}.
\end{align*}
Choose $J_0$ such that $J_0 \geq (2r'pC_\ep)^{\frac{n-m}n}$, by using Adams inequality as above, we then have for any $j \geq J_0$
\[
\int_{\R^n} \Phi_{n,m}(r'p_1C_\ep \beta(n,m) |u_j-\tilde{v}_j|^{\frac n{n-m}})dx \leq C
\]
Using Lemma \ref{expointegral} we obtain
\begin{equation}\label{eq:Omegaj22}
\sup_{j\in \N} \int_{\Om_{j,2}^2} \Phi_{n,m}(p_1 \beta(n,m) |u_j|^{\frac n{n-m}}) dx \leq C\sup_{j\in \N} \int_{\R^n} \Phi_{n,m}(p_2 \beta(n,m)|\tilde{v}_j|^{\frac n{n-m}}) dx.
\end{equation} 
Combining \eqref{eq:mienkobichan}, \eqref{eq:Omegaj1}, \eqref{eq:Omegaj21} and \eqref{eq:Omegaj22} we obtain the conclusion of this lemma. 
\end{proof}

In the proof of Theorem \ref{Maintheoremfull} below, we will need the following technical lemmas. The first one is a radial lemma which asserts that 
\begin{lemma}\label{Radiallemma}
Let $u \in L^{p}(\R^n)$, then 
\[
u^\sharp(x) \leq \lt(\frac 1{|x|^n \om_n}\rt)^{\frac 1p} \lt(\int_{\R^n} |u(y)|^p dy\rt)^{\frac 1p}.
\]
\end{lemma}
The second lemma reads as follows
\begin{lemma}\label{truncation}
Let $u\in L^{\frac nm}(\R^n)$ such that $\|u\|_{\frac nm} \leq 1$, then for any $ p>0$, and $R > 0$ there exists constant $C(n,m,p,R)$ depending only on $n,m,p,$ and $R$ such that 
\[
\int\limits_{|x|\geq R} \Phi_{n,m}(p \beta(n,m) |u^\sharp|^{\frac n{n-m}}) dx \leq C(n,m,p,R).
\]
\end{lemma}
\begin{proof} 
By Radial lemma, we have for any $|x|\geq R$ that 
\[
u^{\sharp}(x) \leq u^\sharp(R) \leq \frac1{|R|^m} \lt(\frac 1{\om_n}\rt)^{\frac mn}.
\]
It is easy to see that for any $A > 0$, the exists a constant $C_{n,m,A}$ such that 
\[
\Phi_{n,m}(t) \leq C_{n,m,A} t^{\frac nm -1},\quad \forall\,t \leq A,
\]
hence there exists $C(n,m,p,R)$ such that 
\[
\Phi_{n,m}(p \beta(n,m) |u^\sharp(x)|^{\frac n{n-m}}) \leq C(n,m,p,R) |u^\sharp(x)|^{\frac nm},\quad \forall\, |x|\geq R.
\]
By this inequality, we obtain the conclusion of Lemma \ref{truncation}.
\end{proof}

We proceed our proof by contradiction argument as in the proof of Theorem \ref{Maintheorem}. Suppose that there exists a sequence $\{u_j\}_j \subset C_0^\infty(\R^n)$ such that $\| u_j\|_{m,\frac nm} \leq 1$, $u_j$ converges weakly to a nonzero function $u$ in $W^{m,\frac nm}(\R^n)$, and a number $p_1 \in (1,Q_{n,m}(u))$ such that 
\begin{equation}\label{eq:contraassumpunboundedfull}
\lim_{j\to\infty} \int_{\R^n} \Phi_{n,m}(p \beta(n,m) |u_j|^{\frac n{n-m}}) dx =\infty.
\end{equation} 
Our aim is to look for a contradiction to \eqref{eq:contraassumpunboundedfull}. Using Rellich--Kondrachov theorem, by passing to a subsequence if necessary, we can assume that $u_j$ converges a.e to $u$ in $\R^n$ and also converges in $L^p_{loc}(\R^n)$ for any $p < \infty$. In the case $m$ odd, we can additionally assume that $(-\Delta+I)^{\frac{m-1}2} u_j $ converges a.e to $(-\Delta+I)^{\frac{m-1}2} u$ in $\R^n$ and also converges in $L^p_{loc}(\R^n)$ for any $p < \frac n{m-1}$.

We write the integral in \eqref{eq:contraassumpunboundedfull} as
\begin{align*}
\int_{\R^n} \Phi_{n,m}(p_1 \beta(n,m) |u_j|^{\frac n{n-m}}) dx & = \int_{\R^n} \Phi_{n,m}(p_1 \beta(n,m) |u_j^\sharp|^{\frac n{n-m}}) dx\\
&=\int_{B_R} \Phi_{n,m}(p_1 \beta(n,m) |u_j^\sharp|^{\frac n{n-m}}) dx \\
&\quad\quad\quad + \int_{B_R^c} \Phi_{n,m}(p_1 \beta(n,m) |u_j^\sharp|^{\frac n{n-m}}) dx,
\end{align*}
where $B_R =\{x\, :\, |x|\leq  R\}$ and $B_R^c = \R^n \setminus B_R$. Lemma \ref{truncation} and our assumption \eqref{eq:contraassumpunboundedfull} imply that
\begin{equation}\label{eq:cutoffintegralfull}
\lim_{j\to \infty} \int_{B_R} \Phi_{n,m}(p_1 \beta(n,m) |u_j^\sharp|^{\frac n{n-m}}) dx =\infty.
\end{equation}
Note that, for $l < \frac nm -1$ we have 
\[
\int_{B_R} (u_j^\sharp)^{l\frac{n}{n-m}} dx \leq \lt(\int_{B_R} (u_j^\sharp)^{\frac nm} dV_g\rt)^{\frac{ml}{n-m}} |B_R|^{1 -\frac{ml}{n-m}} \leq C(n,m,R).
\]
These latter inequalities and \eqref{eq:cutoffintegralfull} imply
\begin{equation}\label{eq:cutoffintegral}
\lim_{j\to \infty} \int_0^{|B_R|} \exp(p_1 \beta(n,m) |u_j^*(s)|^{\frac n{n-m}}) ds = \lim_{j\to \infty} \int_{B_R} \exp(p_1 \beta(n,m) |u_j^\sharp|^{\frac n{n-m}}) dx =\infty.
\end{equation}
We divide our proof into two cases:

{\bf Case $1$: Suppose that $m$ is even.} We can express $m=2k$ for some $k\geq 1$. Denote $f_j = (-\Delta+I)^k u_j$ and $f =(-\Delta+I)^k u$ then $f_j$ converges weakly to $f$ in $L^{\frac n{2k}}(\R^n)$. It follows from Lemma $2$ in \cite{doOMacedo2014} that, by passing to a subsequence, we can assume that $f_j^*$ converges a.e to a function $g$ in $(0,\infty)$ and 

%We argue by contradiction. Suppose that there exists a sequence $\{u_j\}_j \subset C_0^\infty(\R^n)$ with $\|u_j\|_{2k,\frac n{2k}} \leq 1$, $u_j$ converges weakly to $u\not\equiv 0$ in $W^{2k,\frac n{2k}}(\R^n)$ and a number $p_1\in (1,Q_{n,m}(u))$ such that
%\begin{equation}\label{eq:giathietphanchung1}
%\lim_{j\to\infty} \int_{\R^n} \Phi_{n,2k}(p_1 \beta(n,2k) |u_j|^{\frac n{n-2k}}) dx = \infty.
%\end{equation}
%By Rellich-Kondrachov theorem, we can assume that,  by passing to a subsequence if necessary, $u_j$ converges to $u$ in $L_{loc}^{\frac n{2k}}(\R^n)$ and almost everywhere. Denote $f_j = \Delta^k u_j$ and $f =\Delta^k u$ then $f_j$ converges weakly to $f$ in $L^{\frac n{2k}}(\R^n)$. It follows from Lemma $2$ in \cite{doOMacedo2014} that, by passing to a subsequence, we can assume that $f_j^*$ converges a.e to a function $g$ in $(0,\infty)$ and 
\begin{equation}\label{eq:normofgspace}
\int_0^\infty g(s)^{\frac n{2k}} ds \geq \int_0^\infty (f^*(s))^{\frac n{2k}} ds = \int_{\R^n} |f(x)|^{\frac n{2k}} dx.
\end{equation}

%Fix $R > 0$, we write
%\begin{align}\label{eq:tachtichphan1}
%\int_{\R^n} \Phi_{n,2k}(p_1 \beta(n,2k) |u_j|^{\frac n{n-2k}}) dx &= \int_{B_R} \Phi_{n,2k}(p_1 \beta(n,2k) |u_j^\sharp|^{\frac n{n-2k}}) dx\notag\\
%&\quad\quad\quad + \int\limits_{|x|\geq R} \Phi_{n,2k}(p_1 \beta(n,2k) |u_j^\sharp|^{\frac n{n-2k}}) dx.
%\end{align}
%Since for any $j \leq \frac n{2k} -1$, we have by applying H\"older inequality that
%\[
%\int_{B_R} |u_j^\sharp|^{j \frac n{n-2k}} dx \leq \lt(\int_{B_R} |u_j^\sharp|^{\frac n{2k}} dx\rt)^{\frac{2kj}{n-2k}} |B_R|^{1-\frac{2kj}{n-2k}} \leq C_j(R).
%\]
%The latter inequality together Lemma \ref{truncation}, \eqref{eq:giathietphanchung1} and \eqref{eq:tachtichphan1} imply that
%\begin{equation}\label{eq:reducetoball1}
%\lim_{j\to\infty} \int_0^{|B_R|} e^{p_1 \beta(n,2k) |u_j^*(s)|^{\frac n{n-2k}}} ds =\lim_{j\to\infty} \int_{B_R} e^{(p_1 \beta(n,2k) |u_j^\sharp|^{\frac n{n-2k}}} dx =\infty.
%\end{equation}
Arguing as in the proof of Theorem \ref{Maintheorem} for $m$ even by using Proposition \ref{keyfull}, we obtain that 
\[
u_j^*(t_1) -u_j^*(t_2) \leq \frac n{n-2k} \frac{c(n,k)}{(n\om_n^{1/n})^{2k}} \int_{t_1}^{t_2} \frac{f_j^*(s)}{s^{1-\frac{2k}n}} ds + C,
\]
for any $0 < t_1 < t_2 < \infty$, with $C$ depends only on $n,k$. Define the function $v_j$ on $(0,|B_R|]$ by
\[
v_j(t) = \frac n{n-2k} \frac{c(n,k)}{(n\om_n^{1/n})^{2k}} \int_{t}^{|B_R|} \frac{f_j^*(s)}{s^{1-\frac{2k}n}} ds,
\]
we then have $v_j(|B_R|) =0$, and
\[
u_j^*(t) \leq v_j(t) + C + u_j^*(|B_R|) \leq v_j(t) + C(n,k,R),
\]
with $C(n,k,R)$ depends only on $n,k,$ and $R$. Choose $\ep > 0$ small enough such that $\overline{p}_1 = (1+\ep) p_1 < Q(n,2k)(u)$, and applying elementary inequality \eqref{eq:elementary} for $p = \frac{n}{n-2k}$, we obtain
\[
(u_j^*(t))^{\frac n{n-2k}} \leq (1+\ep) v_j(t)^{\frac n{n-2k}} + C_\ep C(n,k,R)^{\frac n{n-2k}},\quad\forall\, t \in (0,|B_R|).
\]
Consequently, we must have from \eqref{eq:cutoffintegral} that
\begin{equation}\label{eq:reducetov1}
\lim_{j\to\infty} \int_0^{|B_R|} e^{\overline{p}_1 \beta(n,2k) |v_j(s)|^{\frac n{n-2k}}} ds = \infty.
\end{equation}
On the other hand, using H\"older inequality, we have for any $s\in (0, |B_R|)$
\begin{equation}\label{eq:pointwisespace1}
v_j(s) \leq \lt(\frac1{\beta(n,2k)} \ln \lt(\frac{|B_R|} s\rt)\rt)^{\frac{n-2k}n}.
\end{equation}
Using the same argument in the proof of Theorem \ref{Maintheorem} (the case $m$ even), we obtain that for each $p_2\in (\op_1,Q_{n,2k}(u))$, for eack $j_0\in \N$ and $s_0 \in (0, |B_R|)$, there exist $j\geq j_0$ and $s\in (0,s_0)$ such that
\begin{equation*}\label{eq:reverseinequality1}
v_j(s) \geq \lt(\frac1{p_2\beta(n,2k)} \ln \lt(\frac{|B_R|} s\rt)\rt)^{\frac{n-2k}n}.
\end{equation*}
Thus, up to a subsequence, we can further assume that there exists $s_j \in (0, |B_R|)$ such that
\begin{equation}\label{eq:reversesequence}
v_j(s_j) \geq \lt(\frac1{p_2\beta(n,2k)} \ln \lt(\frac{|B_R|} {s_j}\rt)\rt)^{\frac{n-2k}n}, \quad\text{ and }\quad s_j \leq \frac 1j.
\end{equation}
For any $L> 0$, considering again the truncation operators $T_L$ and $T^L$ defined as above. We then have $T^L(f_j^*)$ and $T_L(f_j^*)$ converge a.e to $T^L(g)$ and $T_L(g)$ in $(0,\infty)$, respectively.

%Fix $p_3\in (p_2, Q_{n,2k}(u))$, we have two following cases:

%{\bf Case $1$: $0 < \int_{\R^n} |u|^{\frac n{2k}} dx + \int_{\R^n} |f|^{\frac n{2k}} dx <1$.} Since $\lim_{L\to\infty} \int_0^\infty T_L(g)^{\frac n{2k}} ds =0$, then we can choose $L > 0$ large enough such that 
%\begin{equation}\label{eq:Lbigger}
%\frac{1 -\int_{\R^n} |u|^{\frac n{2k}} dx - \int_{\R^n} |f|^{\frac n{2k}} dx}{1 -\int_{\R^n} |u|^{\frac n{2k}} dx -\int_0^\infty(g^{\frac n{2k}} -T_L(g)^{\frac n{2k}})ds} > \lt(\frac{p_3}{Q_{n,2k}(u)}\rt)^{\frac{n-2k}{2k}}.
% \end{equation}
%Fix such a constant $L$. 

By \eqref{eq:reversesequence}, there exists $j_0 >0$ such that $v_j(s_j) > L$ for any $j\geq j_0$, hence we can find $r_j \in (s_j, |B_R|)$ such that $v_j(r_j) = L$. Consequently, we have
\begin{align*}
\lt(\frac1{p_2\beta(n,2k)} \ln \lt(\frac{|B_R|} {s_j}\rt)\rt)^{\frac{n-2k}n} -L &\leq v_j(s_j) -v_j(r_j)\\
& =\frac{n}{n-2k} \frac{c(n,k)}{(n\om_n^{1/n})^{2k}} \int_{s_j}^{r_j} \frac{f_j^*(s)}{s^{1-\frac{2k}n}} ds\\
&\leq \frac{n^2}{2k(n-2k)} \frac{c(n,k)}{(n\om_n^{1/n})^{2k}} f^*_j(s_j) |B_R|^{\frac{2k}n}.
\end{align*}
This shows that $\lim_{j\to\infty} f_j^*(s_j) = \infty$, hence there is $j_1 \geq j_0$ such that $f_j^*(s_j) > L$ for any $j\geq j_1$. Hence there exists $t_j \in (s_j,\infty)$ such that $f_j^*(t_j) = L$ and $f_j(s) < L$ for any $s > t_j$. Let $a_j = \min\{r_j,t_j\}$ we then have
\begin{align*}
\lt(\frac1{p_2\beta(n,2k)} \ln \lt(\frac{|B_R|} {s_j}\rt)\rt)^{\frac{n-2k}n} -L& \leq \frac{n}{n-2k} \frac{c(n,k)}{(n\om_n^{1/n})^{2k}} \int_{s_j}^{a_j} \frac{f_j^*(s)-L}{s^{1-\frac{2k}n}} ds\\
&\quad\quad +\frac{n}{n-2k} \frac{c(n,k)}{(n\om_n^{1/n})^{2k}} \int_{s_j}^{r_j} \frac{L}{s^{1-\frac{2k}n}} ds\\
&\leq \lt(\int_0^\infty |T_L(f_j^*)|^{\frac{n}{2k}} ds\rt)^{\frac n{2k}}\lt(\frac1{\beta(n,2k)} \ln \lt(\frac{|B_R|} {s_j}\rt)\rt)^{\frac{n-2k}n}\\
&\quad\quad  + \frac{c(n,k+1)}{(n\om_n^{1/n})^{2k}} L (|B_R|^{\frac{2k}n} -s_j^{\frac{2k}n}),
\end{align*}
here we apply H\"older inequality. 

For any $p_3 \in (p_2,Q_{n,2k}(u))$, there exists $j_2 \geq j_1$ such that 
\[
p_3^{-\frac{n-2k}{2k}} \leq \int_0^\infty |T_L(f_j^*)|^{\frac{n}{2k}} ds,\quad\forall\, j\geq j_2.
\]
For such $j$, we have
\begin{align*}
1-p_3^{-\frac{n-2k}{2k}} \geq  \int_{0}^\infty ((f_j^*)^{\frac n{2k}} -T_L(f_j^*)^{\frac n{2k}})ds.
\end{align*}
Let $j\to \infty$ and using Fatou lemma, we get
\begin{equation*}\label{eq:done1}
1-p_3^{-\frac{n-2k}{2k}} \geq \int_{0}^\infty (g^{\frac n{2k}} -T_L(g)^{\frac n{2k}})ds.
\end{equation*}
Let $L$ tend to $\infty$, we obtain
\[
1-p_3^{-\frac{n-2k}{2k}} \geq \int_{0}^\infty g^{\frac n{2k}}ds \geq \int_{\R^n} |(-\Delta +I)^k u|^{\frac n{2k}} dx
\]
or equivalently
\[
p_3^{-\frac{n-2k}{2k}} \leq 1-\int_{\R^n} |(-\Delta +I)^k u|^{\frac n{2k}} dx,
\]
for any $p_3\in (p_2,Q_{n,2k}(u))$, which is impossible.

%Combining \eqref{eq:Lbigger} and \eqref{eq:done1} we obtain 
%\begin{align*}
%p_3 & \geq \lt(\frac1{1 -\int_{\R^n} |u|^{\frac n{2k}} dx- \int_{0}^\infty (g^{\frac n{2k}} -T_L(g)^{\frac n{2k}})ds}\rt)^{\frac{2k}{n-2k}}\\
%& > \frac{p_3}{Q_{n,2k}(u)} \lt(\frac1{1 - \int_{\R^n} |u|^{\frac n{2k}} dx- \int_{\R^n} |f|^{\frac n{2k}} dx}\rt)^{\frac{2k}{n-2k}}\\
%& =p_3
%\end{align*}
%which is a contradiction.

%{\bf Case $2$: $\int_{\R^n} |u|^{\frac n{2k}} dx + \int_{\R^n} |f|^{\frac n{2k}} dx =1$.} In this case, we must have 
%\[
%\int_{\R^n} |u|^{\frac n{2k}} dx + \int_0^{\infty} g(s)^{\frac n{2k}} ds =1.
%ù\]
%We hence can choose a $L > 0$ such that 
%\[
%\int_{\R^n} |u|^{\frac n{2k}} dx + \int_0^{\infty} (g^{\frac n{2k}} -T_L(g)^{\frac n{2k}}) ds > 1 -\frac12 \lt(\frac 1{p_3}\rt)^{\frac{n-2k}{2k}}.
%\]
%Fix such a $L$. The argument in {\bf Case $1$} shows that \eqref{eq:done1} also holds, hence we obtain a contradiction that 
%\[
%1-\lt(\frac 1{p_3}\rt)^{\frac{n-2k}{2k}} \geq \int_{\R^n} |u|^{\frac n{2k}} dx+ \int_{0}^\infty (g^{\frac n{2k}} -T_L(g)^{\frac n{2k}})ds > 1 -\frac12 \lt(\frac 1{p_3}\rt)^{\frac{n-2k}{2k}}.
%\]
%This finishes our proof of Theorem \ref{Maintheoremfull} when $m$ is even.
 
{\bf Case $2$: Suppose that $m$ is odd.} In this case we can express $m =2k+1$ for some $k\geq 0$.
Since the case $m=1$ was proved in \cite{doO2014}, hence we will treat only the case $m=2k +1, k\geq 1$. Denote $f_j = (-\Delta+I)^k u_j$ and $f =(-\Delta+I)^k u$. We then have $f_j \in W^{1,\frac n{2k+1}}(\R^n)$, $f_j$ converges weakly to $f$ in $W^{1,\frac n{2k+1}}(\R^n)$ and $f_j$ converges a.e to $f$ in $\R^n$. It implies from Sobolev inequality that $\|f_j\|_{\frac n{2k}} \leq C$ for some constant $C$ depending only on $n,k$. 

%Suppose that there exist a sequence $\{u_j\}_j\subset C_0^\infty(\R^n)$ with $\int_\R^n |u_j|^{\frac n{2k+1}} dx + \int_{\R^n} |\nabla\Delta^k u_j|^{\frac n{2k+1}} dx \leq 1$, $u_j$ converges weakly to $u\not\equiv 0$ in $W_0^{2k+1,\frac n{2k+1}}(\R^n)$ and a number $p_1 \in (1, Q_{n,2k+1}(u))$ such that
%\begin{equation}\label{eq:giathietphanchungle}
%\lim_{j\to \infty} \int_{\R^n} \Phi_{n,2k+1}(p_1 \be(n,2k+1) |u_j|^{\frac n{n-2k-1}}) dx = \infty.
%\end{equation}

Repeating the proof of \eqref{eq:welldone} and using again Lemma \ref{keyfull} imply the existence of a positive constant $C_1$ depending only on $n,k$ such that 
\begin{equation}\label{eq:welldonespace}
u_j^*(t_1) -u_j^*(t_2) \leq \frac{c(n,k+1)}{(n\om_n^{1/n})^{2k}} \int_{t_1}^{t_2} (-f^*_j)'(s) s^{\frac{2k}n}  ds +C_1,
\end{equation}
for any $0 < t_1 < t_2 < \infty$.

%Fix $R > 0$, we write
%\begin{align}\label{eq:tachtichphan2}
%\int_{\R^n} \Phi_{n,2k+1}(p_1 \beta(n,2k+1) |u_j|^{\frac n{n-2k-1}}) dx &= \int_{B_R} \Phi_{n,2k}(p_1 \beta(n,2k+1) |u_j^\sharp|^{\frac n{n-2k-1}}) dx\notag\\
%&\quad\quad\quad + \int\limits_{|x|\geq R} \Phi_{n,2k+1}(p_1 \beta(n,2k+1) |u_j^\sharp|^{\frac n{n-2k-1}}) dx.
%\end{align}
%Since for any $j \leq \frac n{2k+1} -1$, we have by applying H\"older inequality that
%\[
%\int_{B_R} |u_j^\sharp|^{j \frac n{n-2k-1}} dx \leq \lt(\int_{B_R} |u_j^\sharp|^{\frac n{2k+1}} dx\rt)^{\frac{(2k+1)j}{n-2k-1}} |B_R|^{1-\frac{(2k+1)j}{n-2k-1}} \leq C_j(R).
%\]
%The latter inequality together Lemma \ref{truncation}, \eqref{eq:giathietphanchungle} and \eqref{eq:tachtichphan2} imply that
%\begin{equation}\label{eq:reducetoball2}
%\lim_{j\to\infty} \int_0^{|B_R|} e^{p_1 \beta(n,2k+1) |u_j^*(s)|^{\frac n{n-2k-1}}} ds =\lim_{j\to\infty} \int_{B_R} e^{(p_1 \beta(n,2k+1) |u_j^\sharp|^{\frac n{n-2k-1}}} dx =\infty.
%\end{equation}
Define the function $v_j$ on $(0,|B_R|)$ by
\[
v_j(t) = \frac{c(n,k+1)}{(n\om_n^{1/n})^{2k}} \int_{t}^{|B_R|} (-f^*_j)'(s) s^{\frac{2k}n}  ds, \quad 0 < t \leq |B_R|,
\]
then $v_j(|B_R|) =0$, and by Radial lemma, we get $u_j^*(t) \leq v_j(t) + C(n,k,R)$ with $C(n,k,R)$ depends only on $n,k,R$. The same arguments in the proof of \eqref{eq:reducetov1} imply that for any $\op_1 \in (p_1,Q_{n,2k+1}(u))$ we have
%Choose $\ep > 0$ small enough such that $\op_1 =(1+\ep) p_1 < Q_{n,2k+1}(u)$, using the elementary inequality \eqref{eq:elementary} for $p = \frac n{n-2k-1}$, we have
%\[
%(u_j^*(t))^{\frac n{n-2k-1}} \leq (1+\ep) v_j(t)^{\frac n{n-2k-1}} + C_\ep C(n,k,R)^{\frac n{n-2k-1}}.
%\]
%The latter inequality and \eqref{eq:reducetoball2} yields 
\begin{equation}\label{eq:forvfunction}
\lim_{j\to\infty} \int_0^{|B_R|} e^{\op_1 \beta(n,2k+1) |v_j(s)|^{\frac n{n-2k-1}}} ds = \infty.
\end{equation}

%Recall that 
%\[
%\int_0^{\infty} \lt(n\om_n^{\frac1n} (-f_j^*)'(s)\rt)^{\frac n{2k+1}} s^{\frac{n-1}{2k+1}} ds =  \int_{\R^n} |\nabla f_j^\sharp|^{\frac n{2k+1}} dx \leq \int_{\R^n} |\nabla f_j|^{\frac n{2k+1}} dx \leq 1,
%\]
%here we use P\'olya-Szeg\"o principle. 
%Applying H\"older inequality to $v_j$, we obtain for any $s\in (0, |B_R|)$ that
%\begin{equation}\label{eq:boundforvj}
%v_j(s) \leq \lt(\frac 1{\beta(n,2k+1)} \ln\lt(\frac{|B_R|}s\rt)\rt)^{\frac{n-2k-1}{n}}.
%\end{equation}
Fix a $\op_1 \in (p_1,Q_{n,2k+1}(u))$. The same arguments in the proof of Theorem \ref{Maintheorem} shows that for each $p_2\in (\op_1,Q_{n,2k+1}(u))$, for each $j_0\in \N$ and $s_0 \in (0, |B_R|)$, there exist $j\geq j_0$ and $s\in (0,s_0)$ such that
\begin{equation*}\label{eq:reverseinequality2}
v_j(s) \geq \lt(\frac1{p_2\beta(n,2k+1)} \ln \lt(\frac{|B_R|} s\rt)\rt)^{\frac{n-2k-1}n}.
\end{equation*}
Thus, up to a subsequence, we can assume that there exists $s_j \in (0, |B_R|)$ such that
\begin{equation}\label{eq:reversesequence2}
v_j(s_j) \geq \lt(\frac1{p_2\beta(n,2k+1)} \ln \lt(\frac{|B_R|} {s_j}\rt)\rt)^{\frac{n-2k-1}n}, \quad\text{ and }\quad s_j \leq \frac 1j.
\end{equation}

For any $L>0$, let us reconsider the truncation operators $T_L$ and $T^L$ defined above, we then have $T^L(f_j)$ and $T_L(f_j)$ converge a.e to $T^L(f)$ and $T_L(f)$ in $\R^n$, respectively. We also have $T^L(f_j)$ and $T_L(f_j)$ converge weakly to $T^L(f)$ and $T_L(f)$ in $W^{1,\frac n{
k+}}(\R^n)$, respectively (see the proof of the case $m$ odd in Theorem \ref{Maintheorem}). 

It follows from \eqref{eq:reversesequence2} that there exists $j_0$ such that $v_j(s_j) > L$ for any $j\geq j_0$, hence there exists $r_j \in (s_j, |B_R|)$ such that $v_j(r_j) = L$. On the other hand, we have
\[
v_j(s_j) \leq \frac{c(n,k+1) |B_R|^{\frac{2k}n}}{(n\om_n^{1/n})^{2k}} (f_j^*(s_j) -f_j^*(B_R)).
\]
Then we must have $\lim_{j\to \infty} f_j^*(s_j) = \infty$ by \eqref{eq:reversesequence2}. So there exists $j_1\geq j_0$ such that $f_j^*(s_j) >L$ for any $j\geq j_1$. For such $j$, there is $t_j \in (s_j,\infty)$ such that $f_j^*(t_j) = L$ and $f^*_j(s)< L$ for any $s > t_j$. Hence
\[
\int_{\R^n} |\nabla T_L(f_j^\sharp)|^{\frac n{2k+1}} dx = \int_0^{t_j} \lt(n\om_n^{\frac1n} (-f_j^*)'(s)\rt)^{\frac n{2k+1}} s^{\frac{n-1}{2k+1}} ds.
\]
Let $a_j =\min\{t_j, r_j\}$, we then have
\begin{align}\label{eq:hoangthom}
v_j(s_j) -v_j(r_j) &= \frac{c(n,k+1)}{(n\om_n^{1/n})^{2k}}\int_{s_j}^{r_j} (-f_j^*)'(s) s^{\frac{2k}n} ds\notag\\
&=\frac{c(n,k+1)}{(n\om_n^{1/n})^{2k}}\int_{s_j}^{a_j} (-f_j^*)'(s) s^{\frac{2k}n} ds +\frac{c(n,k+1)}{(n\om_n^{1/n})^{2k}}\int_{a_j}^{r_j} (-f_j^*)'(s) s^{\frac{2k}n} ds\notag\\
&\leq \lt(\int_0^{t_j} \lt(n\om_n^{\frac1n} (-f_j^*)'(s)\rt)^{\frac n{2k+1}} s^{\frac{n-1}{2k+1}} ds\rt)^{\frac{2k+1}n} \lt(\frac1{\beta(n,2k+1)} \ln \lt(\frac{a_j}{s_j}\rt)\rt)^{\frac{n-2k-1}n}\notag\\
&\quad\quad \quad  + \frac{c(n,k+1) r_j^{\frac{2k}n}}{(n\om_n^{1/n})^{2k}}(f_j^*(a_j) -f_j^*(r_j))\notag\\
&\leq \lt(\int_{\R^n} |\nabla T_L(f^\sharp)|^{\frac n{2k+1}} dx\rt)^{\frac{2k+1}n} \lt(\frac1{\beta(n,2k+1)} \ln \lt(\frac{|B_R|}{s_j}\rt)\rt)^{\frac{n-2k-1}n}\notag\\
&\quad\quad\quad + \frac{c(n,k+1) |B_R|^{\frac{2k}n}}{(n\om_n^{1/n})^{2k}} L,
\end{align}
here we use the facts that if $t_j< r_j$ then $f_j^*(a_j) -f_j^*(r_j) \leq L$ while if $t_j \geq r_j$ then $f_j^*(a_j) -f_j^*(r_j) =0$. 

For any $p_3 \in (p_2,Q_{n,2k+1}(u))$, it follows from \eqref{eq:reversesequence2}, \eqref{eq:hoangthom} and the fact $v_j(r_j) =L$ for any $j\geq j_1$ that there exists $j_2\geq j_1$ such that
\[
p_3^{-\frac{n-2k-1}{2k+1}} \leq \int_{\R^n} |\nabla T_L(f_j^\sharp)|^{\frac n{2k+1}} dx,\quad\forall\, j\geq j_2.
\]
Note that $(T_L(f_j))^\sharp = T_L(f_j^\sharp)$, hence by P\'olya-Szeg\"o principle, we have
\[
p_3^{-\frac{n-2k-1}{2k+1}} \leq  \int_{\R^n} |\nabla T_L(f_j)|^{\frac n{2k+1}} dx,
\]
or equivalently
\begin{align*}
1-p_3^{-\frac{n-2k-1}{2k+1}}&\geq 1-\int_{\R^n} |\nabla T_L(f_j)|^{\frac n{2k+1}} dx\\
&\geq \int_{\R^n} |f_j|^{\frac{n}{2k+1}} dx +\int_{\R^n} |\nabla f_j|^{\frac n{2k+1}} dx -\int_{\R^n} |\nabla T_L(f_j)|^{\frac n{2k+1}} dx\\
&= \int_{\R^n} |f_j|^{\frac{n}{2k+1}} dx +\int_{\R^n} |\nabla T^L(f_j)|^{\frac n{2k+1}} dx,
\end{align*}
for any $j\geq j_2$. Let $j$ tend to infinity, using Fatou lemma and the weak semicontinuity of the $L^{\frac n{2k+1}}-$ norm of gradient, we then have
\[
1-p_3^{-\frac{n-2k-1}{2k+1}} \geq \int_{\R^n} |f|^{\frac{n}{2k+1}} dx +\int_{\R^n} |\nabla T^L(f)|^{\frac n{2k+1}} dx,
\]
for any $p_3\in (p_2,Q_{n,2k+1}(u))$ and $L>0$. Note that
\[
\int_{\R^n} |\nabla T^L(f)|^{\frac n{2k+1}} dx =\int_{\{|f| < L\}} |\nabla f|^{\frac n{2k+1}} dx,
\]
which then implies
\[
\lim_{L\to \infty} \int_{\R^n} |\nabla T^L(f)|^{\frac n{2k+1}} dx =\int_{\R^n} |\nabla f|^{\frac n{2k+1}} dx.
\]
Let $L$ tend to infinity we obtain
\[
1-p_3^{-\frac{n-2k-1}{2k+1}} \geq \int_{\R^n} |f|^{\frac{n}{2k+1}} dx +\int_{\R^n} |\nabla f|^{\frac n{2k+1}} dx,
\]
or equivalently
\[
p_3^{-\frac{n-2k-1}{2k+1}}\leq 1 - \|u\|_{2k+1,\frac{n}{2k+1}}^{\frac n{2k+1}},
\]
for any $p_3\in (p_2,Q_{n,2k+1}(u))$ which is impossible.

\subsection{Proof of Theorem \ref{Maintheorem2}}
The proof of Theorem \ref{Maintheorem2} follows the same lines in the proof of Theorem \ref{Maintheoremfull}. In the following, we will sketch this proof. It is enough to prove Theorem \ref{Maintheorem2} for compactly supported smooth functions. We argue by contradiction argument. Suppose there exist $\{u_j\}_j\subset C_0^\infty(\R^n)$ such that $\|u_j\| \leq 1$, $u_j$ converges weakly to a nonzero function $u$ in $W^{m,\frac nm}(\R^n)$ and a number $p_1\in (1,R_{n,m}(u))$ such that 
\[
\lim_{j\to\infty} \int_{\R^n} \Phi_{n,m}(p_1\beta(n,m) |u_j|^{\frac n{n-m}}) dx =\infty.
\]
By the same argument in proof of Theorem \ref{Maintheoremfull}, we have for any $R>0$ 
\[
\lim_{j\to\infty} \int_0^{|B_R|} \Phi_{n,m}(p_1 \beta(n,m) |u_j^*(s)|^{\frac n{n-m}})ds =\infty.
\]
We divide our proof into two cases.

{\bf Case $1$: Suppose that $m$ is even.} In this case, we can write $m =2k$ for some $k\geq 1$. Denote $f_j = \Delta^k u_j$ and $f = \Delta^k u$, then $f_j$ converges weakly to $f$ in $L^{\frac n{2k}}(\R^n)$. Using Rellich-Kondrachov theorem, we can assume that $u_j$ converges a.e to $u$ in $\R^n$ and also converges in $L^p_{loc}(\R^n)$ for any $p < \infty$. Possibly passing to a subsequence, we can assume further that $f_j^*$ converges a.e to a function $g$ in $(0,\infty)$ with
\[
\int_0^\infty g(s)^{\frac n{2k}} ds \geq \int_0^\infty (f^*(s))^{\frac n{2k}} ds = \int_{\R^n} |f(x)|^{\frac n{2k}} dx.
\]

We claim that 
\begin{equation}\label{eq:claimnorm}
\int_{\R^n} |u|^{\frac n{2k}} dx + \int_0^\infty g(s)^{\frac n{2k}} ds \leq 1.
\end{equation}
Indeed, for any $R > 0$ we have
\[
1 \geq \int_{B_R} |u_j|^{\frac n{2k}} dx + \int_0^{|B_R|} (f_j^*(s))^{\frac n{2k}} ds.
\]
Since $u_j\to u$ in $L_{loc}^{\frac n{2k}}(\R^n)$ and $f_j^* \to g$ a.e in $(0,\infty)$, then  letting $j\to \infty$ and using Fatou lemma, we obtain
\[
1 \geq \int_{B_R} |u|^{\frac n{2k}} dx + \int_0^{|B_R|} g(s)^{\frac n{2k}} ds.
\]
Letting $R\to \infty$ proves our claim \eqref{eq:claimnorm}.

Define the function $v_j$ on $(0,|B_R|)$ by
\[
v_j(t) = \frac n{n-2k} \frac{c(n,k)}{(n\om_n^{1/n})^{2k}} \int_{t}^{|B_R|} \frac{f_j^*(s)}{s^{1-\frac{2k}n}} ds,
\]
then $v_j(|B_R|) =0$. Repeating the argument in the proof of Theorem \ref{Maintheoremfull} for $m$ even, we have that for any $\op_1 \in (p_1, R_{n,2k}(u))$
\[
\lim_{j\to\infty} \int_0^{|B_R|} e^{\op_1 \beta(n,2k) v_j(s)^{\frac n{n-2k}}} ds = \infty,
\]
(with the helps of Proposition \ref{keypropo}, Radial lemma and the elementary inequality \eqref{eq:elementary}) and that for any $L > 0$ and $p_3 \in (p_2,R_{n,2k}(u))$ with some $p_2 \in (\op_1,R_{n,2k}(u))$, there exists $j_0$ depending only on $p_2, p_3$ and $L$ such that 
\[
p_3^{-\frac{n-2k}{2k}} \leq \int_0^\infty |T_L(f_j^*)|^{\frac{n}{2k}} ds,\quad\forall\, j\geq j_2.
\]
For such $j$, we have
\begin{align*}
1-p_3^{-\frac{n-2k}{2k}} \geq \int_{\R^n} |u_j|^{\frac n{2k}} dx + \int_{0}^\infty ((f_j^*)^{\frac n{2k}} -T_L(f_j^*)^{\frac n{2k}})ds.
\end{align*}
Let $j\to \infty$ and using Fatou lemma, we get
\begin{equation*}
1-p_3^{-\frac{n-2k}{2k}} \geq \int_{\R^n} |u|^{\frac n{2k}} dx+ \int_{0}^\infty (g^{\frac n{2k}} -T_L(g)^{\frac n{2k}})ds.
\end{equation*}
Let $L$ tend to $\infty$, we obtain
\[
1-p_3^{-\frac{n-2k}{2k}} \geq \int_{\R^n} |u|^{\frac n{2k}} dx + \int_{0}^\infty g^{\frac n{2k}}ds \geq \int_{\R^n} |u|^{\frac n{2k}} dx+ \int_{\R^n} |\Delta^k u|^{\frac n{2k}} dx
\]
or equivalently
\[
p_3^{-\frac{n-2k}{2k}} \leq 1-\|u\|^{\frac nm}
\]
for any $p_3\in (p_2,R_{n,2k}(u))$, which is impossible.

{\bf Case $2$: Suppose that $m$ is odd.} In this case, we can write $m=2k+1$ for some $k\geq 0$. We only have to prove for $k >1$. Denote $f_j = \Delta^k u_j$ and $f =\Delta^k u$. We then have $f_j$ converges weakly to $f$ in $W^{1,\frac n{2k+1}}(\R^n)$ and $\|f_j\|_{\frac n{2k}} \leq C$ for some constant $C$ depending only on $n$ and $k$ by the Sobolev inequality. Using Rellich-Kondrachov theorem, we can assume further that $u_j$ and $f_j$ converge a.e to $u$ and $f$ in $\R^n$, respectively (also that $f_j$ converges to $f$ in $L^p_{loc}(\R^n)$ for any $p < \frac n{2k}$).

Define the function $v_j$ on $(0, |B_R|)$ by
\[
v_j(t) = \frac{c(n,k+1)}{(n\om_n^{1/n})^{2k}} \int_{t}^{|B_R|} (-f^*_j)'(s) s^{\frac{2k}n}  ds, \quad 0 < t \leq |B_R|,
\]
the $v_j(|B_R|) =0$. Repeating the argument in the proof of Theorem \ref{Maintheoremfull} for $m$ odd, we have that for any $\op_1 \in (p_1, R_{n,2k+1}(u))$
\[
\lim_{j\to\infty} \int_0^{|B_R|} e^{\op_1 \beta(n,2k+1) v_j(s)^{\frac n{n-2k-1}}} ds = \infty,
\]
(with the helps of Proposition \ref{keypropo}, Radial lemma and the elementary inequality \eqref{eq:elementary}) and that for any $L > 0$ and $p_3 \in (p_2,R_{n,2k+1}(u))$ with some $p_2 \in (\op_1,R_{n,2k+1}(u))$, there exists $j_0$ depending only on $p_2, p_3$ and $L$ such that 
\[
p_3^{-\frac{n-2k-1}{2k+1}} \leq \int_{\R^n} |\nabla T_L(f_j^\sharp)|^{\frac n{2k+1}} dx,\quad\forall\, j\geq j_2.
\]
%Note that $(T_L(f_j))^\sharp = T_L(f_j^\sharp)$, hence by P\'olya-Szeg\"o principle, we have
%\[
%p_3^{-\frac{n-2k-1}{2k+1}} \leq  \int_{\R^n} |\nabla T_L(f_j)|^{\frac n{2k+1}} dx,
%\]
%or equivalently
Then by P\'olya-Szeg\"o principle, we have
\begin{align*}
1-p_3^{-\frac{n-2k-1}{2k+1}}&\geq 1-\int_{\R^n} |\nabla T_L(f_j)|^{\frac n{2k+1}} dx\\
&\geq \int_{\R^n} |u_j|^{\frac{n}{2k+1}} dx +\int_{\R^n} |\nabla f_j|^{\frac n{2k+1}} dx -\int_{\R^n} |\nabla T_L(f_j)|^{\frac n{2k+1}} dx\\
&= \int_{\R^n} |u_j|^{\frac{n}{2k+1}} dx +\int_{\R^n} |\nabla T^L(f_j)|^{\frac n{2k+1}} dx,
\end{align*}
for any $j\geq j_2$. Note that $T^L(f_j)$ converges weakly to $T^L(f)$ in $W^{1,\frac{n}{2k+1}}(\R^n)$ (using the same argument in the proof of Theorem \ref{Maintheorem} in the case $m$ odd).  Let $j$ tend to infinity, using Fatou lemma and the weak semicontinuity of the $L^{\frac n{2k+1}}-$ norm of gradient, we then have
\[
1-p_3^{-\frac{n-2k-1}{2k+1}} \geq \int_{\R^n} |f|^{\frac{n}{2k+1}} dx +\int_{\R^n} |\nabla T^L(f)|^{\frac n{2k+1}} dx,
\]
for any $p_3\in (p_2,R_{n,2k+1}(u))$ and $L>0$.
% Note that
%\[
%\int_{\R^n} |\nabla T^L(f)|^{\frac n{2k+1}} dx =\int_{\{|f| < L\}} |\nabla f|^{\frac n{2k+1}} dx,
%\]
%which then implies
%\[
%\lim_{L\to \infty} \int_{\R^n} |\nabla T^L(f)|^{\frac n{2k+1}} dx =\int_{\R^n} |\nabla f|^{\frac n{2k+1}} dx.
%\]
Let $L$ tend to infinity we obtain
\[
1-p_3^{-\frac{n-2k-1}{2k+1}} \geq \int_{\R^n} |f|^{\frac{n}{2k+1}} dx +\int_{\R^n} |\nabla f|^{\frac n{2k+1}} dx,
\]
or equivalently
\[
p_3^{-\frac{n-2k-1}{2k+1}}\leq 1 - \|u\|^{\frac n{2k+1}},
\]
for any $p_3\in (p_2,R_{n,2k+1}(u))$ which is impossible.

%\subsection{Proof of Corollary \ref{outsidezerofunction}}
%Corollary \ref{outsidezerofunction} is immediate consequence of Theorem \ref{Maintheoremfull} and \ref{Maintheorem2}. We prove it by contradiction argument. Let us first prove \eqref{eq:forTh2}. Suppose that there exist a subsequence $\{u_j\}_j \subset W^{m,\frac nm}(\R^n)\setminus \mW$ such that $\|u_j\|_{m,\frac nm} \leq 1$ and a number $\al\in (0,1)$ such that
%\begin{equation}\label{eq:forTh2contra}
%\lim_{j\to\infty} \int_{\R^n} \Phi_{n,m}((1+\al\|u_j\|_{\frac nm}^{\frac nm})^{\frac1{n-m}} |u_j|^{\frac n{n-m}}) dx =\infty.
%\end{equation}
%By passing to a subsequence, we can assume that $u_j$ converges weakly to a nonzero function $u$ and $u_j$ converges a.e to $u$ in $\R^n$. Fatous lemma implies that 
%\[
%\|u\|_{\frac nm} \leq \liminf_{j\to\infty} \|u_j\|_{\frac nm}
%\]
%Fix $\beta \in (\al,1)$, there exists $j_0$ such that 
%\[
%1+\al\|u\|_{\frac nm}^{\frac nm} \leq 1+ \beta\|u_j\|_{\frac nm}^{\frac nm},\quad\forall\, j\geq j_0.
%\]

\section*{Acknownledgments}
This work is supported by CIMI postdoctoral research fellowship.

\end{document}